\documentclass[12pt]{article}
\usepackage{amssymb,amsthm,amsmath,amsfonts,latexsym,tikz,hyperref,ytableau}
\usepackage[hmargin=1in,vmargin=1in]{geometry}

\newtheorem{thm}{Theorem}[section]
\newtheorem{prop}[thm]{Proposition}
\newtheorem{cor}[thm]{Corollary}
\newtheorem{lem}[thm]{Lemma}
\newtheorem{conj}[thm]{Conjecture}
\newtheorem{exa}[thm]{Example}
\newtheorem{ques}[thm]{Question}

\newcommand{\peq}{\preceq}

\newcommand{\tles}{\triangleleft}

\newcommand{\tleq}{\trianglelefteq}

\DeclareMathOperator{\SSYT}{SSYT}
\DeclareMathOperator{\SSYTS}{SSYTS}
\DeclareMathOperator{\Pk}{Pk}
\DeclareMathOperator{\RGF}{RGF}
\DeclareMathOperator{\NC}{NC}
\DeclareMathOperator{\asm}{asm}
\DeclareMathOperator{\prf}{pf}
\DeclareMathOperator{\Pin}{Pin}
\DeclareMathOperator{\pin}{pin}

\newcommand{\ben}{\begin{enumerate}}
\newcommand{\een}{\end{enumerate}}
\newcommand{\ble}{\begin{lem}}
\newcommand{\ele}{\end{lem}}
\newcommand{\bth}{\begin{thm}}
\renewcommand{\eth}{\end{thm}}
\newcommand{\bpr}{\begin{prop}}
\newcommand{\epr}{\end{prop}}
\newcommand{\bco}{\begin{cor}}
\newcommand{\eco}{\end{cor}}
\newcommand{\bcon}{\begin{conj}}
\newcommand{\econ}{\end{conj}}
\newcommand{\bde}{\begin{defn}}
\newcommand{\ede}{\end{defn}}
\newcommand{\bex}{\begin{exa}}
\newcommand{\eex}{\end{exa}}
\newcommand{\barr}{\begin{array}}
\newcommand{\earr}{\end{array}}
\newcommand{\btab}{\begin{tabular}}
\newcommand{\etab}{\end{tabular}}
\newcommand{\beq}{\begin{equation}}
\newcommand{\eeq}{\end{equation}}
\newcommand{\bea}{\begin{eqnarray*}}
\newcommand{\eea}{\end{eqnarray*}}
\newcommand{\bal}{\begin{align*}}
\newcommand{\bce}{\begin{center}}
\newcommand{\ece}{\end{center}}
\newcommand{\bpi}{\begin{picture}}
\newcommand{\epi}{\end{picture}}
\newcommand{\bpp}{\begin{picture}}
\newcommand{\epp}{\end{picture}}
\newcommand{\bfi}{\begin{figure} \begin{center}}
\newcommand{\efi}{\end{center} \end{figure}}
\newcommand{\bprf}{\begin{proof}}
\newcommand{\eprf}{\end{proof}\medskip}
\newcommand{\capt}{\caption}
\newcommand{\bsl}{\begin{slide}{}}
\newcommand{\esl}{\end{slide}}
\newcommand{\bfr}{\begin{frame}}
\newcommand{\efr}{\end{frame}}

\newcommand{\hqed}{\hfill \qed}

\newcommand{\eqqed}[1]{$\rule{1ex}{0ex}\hfill{\dil#1}\hfill\qed$}

\newcommand{\ol}{\overline}

\newcommand{\hs}[1]{\hspace{#1}}
\newcommand{\hso}[1]{\hspace{-1pt}}
\newcommand{\vs}[1]{\vspace{#1}}

\newcommand{\emp}{\emptyset}

\newcommand{\sbe}{\subseteq}

\newcommand{\setm}{\setminus}

\newcommand{\iso}{\cong}

\newcommand{\ptn}{\vdash}

\newcommand{\jn}{\vee}

\newcommand{\mt}{\wedge}

\newcommand{\case}[4]{\left\{\barr{ll}#1&\mbox{#2}\\#3&\mbox{#4}\earr\right.}

\def\<{\langle}
\def\>{\rangle}

\newcommand{\spn}[1]{\langle{#1}\rangle}

\newcommand{\ra}{\rightarrow}

\newcommand{\be}{\beta}

\newcommand{\io}{\iota}
\newcommand{\ka}{\kappa}
\newcommand{\la}{\lambda}

\newcommand{\si}{\sigma}

\newcommand{\Om}{\Omega}

\newcommand{\1}{{\bf 1}}

\newcommand{\bx}{{\bf x}}

\newcommand{\bbN}{{\mathbb N}}

\newcommand{\bbR}{{\mathbb R}}

\newcommand{\bbZ}{{\mathbb Z}}

\newcommand{\cC}{{\cal C}}
\newcommand{\cD}{{\cal D}}

\newcommand{\cI}{{\cal I}}

\newcommand{\cJ}{{\cal J}}
\newcommand{\cK}{{\cal K}}

\newcommand{\cM}{{\cal M}}

\newcommand{\cO}{{\cal O}}

\newcommand{\cS}{{\cal S}}

\newcommand{\fS}{{\mathfrak S}}

\newcommand{\nb}{\ol{n}}

\DeclareMathOperator{\Av}{Av}

\DeclareMathOperator{\Des}{Des}

\DeclareMathOperator{\Inv}{Inv}

\newcommand{\dil}{\displaystyle}

\begin{document}
\pagestyle{plain}

\title{Log-concavity and log-convexity  via distributive lattices
}
\author{Jinting Liang\\[--5pt]
\small Department of Mathematics, University of British Columbia,\\[-5pt]
\small 1984 Mathematics Rd, Vancouver, BC V6T 1Z2, CANADA, {\tt liangj@math.ubc.ca}\\
and\\
Bruce E. Sagan (corresponding author)
\\[-5pt]
\small Department of Mathematics, Michigan State University,\\[-5pt]
\small East Lansing, MI 48824-1027, USA, {\tt bsagan@msu.edu}
}

\date{\today\\[10pt]
	\begin{flushleft}
	                                       \\[5pt]
	\end{flushleft}}

\maketitle

\begin{abstract}

The FKG inequality is a powerful tool for proving inequalities in distributive lattices.  We show how a special case, which we call the Order Ideal Lemma,  can be used to demonstrate a wide array of log-concavity and log-convexity results in a combinatorial manner.
We use the Order Ideal Lemma to prove log-concavity and log-convexity of various sequences involving lattice paths (Catalan, Motzkin and large Schr\"oder numbers), intervals in Young's lattice, order polynomials, specializations of Schur and Schur $Q$-functions, Lucas sequences, descent and peak polynomials of permutations, pattern avoidance, set partitions, and noncrossing partitions.  We end with a section with conjectures and outlining future directions.

\end{abstract}

\section{Introduction}

Let 
$$
(a_n)_{n\ge0} = a_0, a_1, a_2, \ldots
$$
be a sequence of real numbers.  The sequence is {\em log-concave} if
\beq
\label{lcn}
a_n^2 \ge a_{n-1} a_{n+1}
\eeq
for all $n\ge1$.  A {\em log-convex} sequence is one satisfying
\beq
\label{lcx}
a_n^2\le a_{n-1} a_{n+1}
\eeq
for all $n\ge1$.  Log-concave and log-convex sequences abound in combinatorics, algebra, and geometry.  See the survey articles of Stanley~\cite{sta:lcu}, Brenti~\cite{bre:lcu}, and Br\"and\'en~\cite{bra:ulc} for more information.  By specialization of the FKG inequality~\cite{FKG:cip}, we obtain a combinatorial  tool for proving log-concavity and log-convexity using order ideals in distributive lattices.  We then apply this method to obtain a host of such results both old and new in a unified way.

Let us review some basic concepts from the theory of partially ordered sets (posets).  More comprehensive treatments can be found in the books of Sagan~\cite{sag:aoc} or Stanley~\cite{sta:ec1}.
All of our posets will be finite.
A {\em lower order ideal} in a poset $(P,\preceq)$ is $I\sbe P$ such that   if $x\in I$ and $y\preceq x$ then $y\in I$.  Similarly, an {\em upper order ideal} is $J\sbe P$ satisfying $x\in J$ and $y\succeq x$ implies $y\in J$.
We will use ``order ideal" to refer to a subset which could be either.  
Say that poset $L$ is a {\em lattice} if every pair $x,y\in L$ has a greatest lower bound or {\em meet}, $x\mt y$, as well as a least upper bound or {\em join}, $x\jn y$.  The lattice is {\em distributive} if it satisfies either of the two equivalent distributive laws that, for all $x,y,z\in L$,
$$
x\mt(y\jn z) = (x\mt y) \jn (x\mt z)
$$
and
$$
x\jn(y\mt z) = (x\jn y) \mt (x\jn z).
$$

We can now state our fundamental result which we call the Order Ideal Lemma.  It is an easy consequence of the FKG Inequality~\cite{FKG:cip}.  But since its proof involves concepts not needed in the rest of this paper, we postpone the demonstration until Subsection~\ref{poi}.  For a set $S$ we will use both $|S|$ and $\#S$ for its cardinality.
\ble[The Order  Ideal Lemma]
\label{oil}
Let $L$ be a distributive lattice and suppose that $I,J\sbe L$ are ideals.
\ben
\item[(a)]  If $I,J$ are both lower ideals or both upper ideals then
$$
|I| \cdot |J| \le |I\cap J| \cdot |L|.
$$
\item[(b)]
If one of $I,J$ is a lower order ideal and the other is upper then
$$
|I| \cdot |J| \ge |I\cap J| \cdot |L|.
$$
\een
\ele

We note that part (a) of the Order Ideal Lemma appeared in a paper of Daykin, Kleitman, and West~\cite[equation (2)]{DKW:nmt}.  But it's wide applicability does not seem to have been previously recognized.

Our general strategy for proving log-convexity of a sequence $(a_n)_{n\ge0}$ will be to construct lattices $L_n$ with $|L_n| = a_n$.  If we can find inside $L_{n+1}$ two lower order ideals $I,J$ such 
that $|I|=|J| = a_n$ and $|I\cap J| = a_{n-1}$ then we will be done by part (a) of the Order Ideal Lemma.  Similarly, part (b) can be used to prove log-concavity.

The rest of this paper is structured as follows.  In the next section we will use lattice paths to prove log-convexity of sequences involving the Catalan, Motzkin, and large Schr\"oder numbers.
In Section~\ref{yl} we use various intervals in Young's lattice to give   log-concavity and log-convexity results.  Some of these specialize to the fact that  show that various sequences of binomial coefficients are  log-concave.
We begin Section~\ref{op} by showing that for any poset, the sequence obtained by evaluating its (enriched) order polynomial at nonnegative integers is always log-concave.  As a consequence, we obtain 
log-concavity of sequences of specializations of Schur and Schur $Q$-functions.
Section~\ref{gls} is dedicated to generalized Lucas sequences which are those satisfying $a_n=a_{n-1}+a_{n-2}$ for $n\ge2$.  We show that any such sequence which has positive initial conditions alternates between satisfying~\eqref{lcn} and~\eqref{lcx}.   In Section~\ref{per} we consider various  sequences related to permutations, including those defined by pattern avoidance as well as sequences of descent and peak polynomials.  The focus of Section~\ref{sp} is set partitions and we show log-concavity of sequences involving Stirling numbers of the second  kind and Narayana numbers.
The last section contains a proof of the Order Ideal Lemma as well as directions for future research.

\section{Lattice paths}
\label{lp}

\bfi
\bce
\begin{tikzpicture}[scale=.5]
\draw (0,1)--(-2,2) (6,1)--(8,2) (-3,4)--(0,5) (9,4)--(6,5) (3,7.5)--(3,8.5);
\filldraw (0,0)--(1,1)--(2,0)--(3,1)--(4,0)--(5,1)--(6,0)--(0,0);
\draw(3,-1) node{$UDUDUD$};
\begin{scope}[shift={(-7,3)}]
\filldraw (0,0)--(1,1)--(2,2)--(3,1)--(4,0)--(5,1)--(6,0)--(0,0);
\draw(3,-1) node{$UUDDUD$};
\end{scope}
\begin{scope}[shift={(7,3)}]
\filldraw (0,0)--(1,1)--(2,0)--(3,1)--(4,2)--(5,1)--(6,0)--(0,0);
\draw(3,-1) node{$UDUUDD$};
\end{scope}
\begin{scope}[shift={(0,6)}]
\filldraw (0,0)--(1,1)--(2,2)--(3,1)--(4,2)--(5,1)--(6,0)--(0,0);
\draw(3,-1) node{$UUDUDD$};
\end{scope}
\begin{scope}[shift={(0,10)}]
\filldraw (0,0)--(1,1)--(2,2)--(3,3)--(4,2)--(5,1)--(6,0)--(0,0);
\draw(3,-1) node{$UUUDDD$};
\end{scope}
\end{tikzpicture}
\ece
\capt{ The poset $\cD_3$}
\label{D3} 
\efi

In this section we will use lattice paths together with the Order Ideal Lemma to give unified proofs of the log-convexity of the sequences of Catalan, Motzkin, and large Schr\"oder numbers.  We begin with a review of some basic definitions.

A {\em lattice path} is a sequence $P:p_0,p_1,\ldots,p_n$ of points in the integer lattice so $p_i\in\bbZ^2$ for all $i$.  A {\em step} of $P$ is the vector $[x_i,y_i]$ from $p_{i-1}$ to $p_i$.  
When the initial vertex $p_0$ of $P$ is known, we can specify $P$ by listing its steps.  
An {\em up step} is a step $U=[1,1]$ and a {\em down step} is $D=[1,-1]$.  A {\em Dyck path of semilength $n$} is a lattice path $P$ satisfying
\ben
\item $P$ starts at $p_0=(0,0)$ and ends at $p_{2n}=(2n,0)$,
\item $P$ uses steps $U$ and $D$ and never goes below the $x$-axis.
\een
Figure~\ref{D3} displays the five Dyck paths of semilength $3$.  Let
$$
\cD_n = \{P \mid \text{$P$ is a Dyck path of semilength $n$}\}.
$$
It is well known that the cardinality of $\cD_n$ is the Catalan number
\beq
\label{Cn}
C_n = \frac{1}{n+1}\binom{2n}{n}.
\eeq

We wish to turn $\cD_n$ into a distributive lattice.  If $P\in\cD_n$ then we let $A(P)$ be the {\em area} of $P$ which is the set of all points of $\bbR^2$ between $P$ and the $x$-axis.  In Figure~\ref{D3} the areas are shaded.  We now define a partial order on $\cD_n$ by
\beq
\label{AP}
P\preceq Q \text{ if and only if } A(P) \subseteq A(Q).
\eeq
The Hasse diagram for $\cD_3$ is in Figure~\ref{D3}.  We note that $\cD_n$ is a distributive lattice.
This follows from the fact that it is isomorphic to an interval in Young's lattice and we will discuss various interesting intervals in the next section.
It is also a consequence of a more general theorem of Ferrari and Pinzani~\cite{FP:llp} giving a criterion for a family of lattice paths ordered by~\eqref{AP} to be a distributive lattice.
The next result follows from easy algebraic manipulations of~\eqref{Cn}.  But our proof is combinatorial and will generalize to other families of paths where a closed-form formula is not known.
We note that a $q$-analogue of this result has been proved by Butler and Flanigan~\cite{BF:nlc}, although it does not use the standard definition of $q$-log convexity discussed in Subsection~\ref{qa} below.
\bth
\label{CnThm}
The sequence $(C_n)_{n\ge0}$ of Catalan numbers is log-convex.
\eth
\bprf
We begin with the distributive lattice $\cD_{n+1}$ and note that $|\cD_{n+1}|= C_{n+1}$.

Let
$$
I=\{P\in\cD_{n+1} \mid \text{$P=UDP'$ for some translated Dyck path $P'$ of semilength $n$}\}.
$$
it follows that $I$ is a lower order ideal because if $P\in I$ and $Q\preceq P$ then~\eqref{AP} forces $Q=UDQ'$ for some $Q'$.
Furthermore, we have an isomorphism of posets $I\iso\cD_n$ given by $P=UDP'\mapsto P'$.  Thus $|I|=C_n$.

Now consider
$$
J=\{P\in\cD_{n+1} \mid \text{$P=P' UD$ for some  Dyck path $P'$ of semilength $n$}\}.
$$
Similar considerations to those in the previous paragraph show that $|J|=C_n$.
Furthermore
$$
I\cap J = \{P\in\cD_{n+1} \mid \text{$P=UD P' UD$ for some  Dyck path $P'$ of semilength $n-1$}\}
$$
so that $|I\cap J| = C_{n-1}$.  Now applying part (a) of the Order Ideal Lemma gives
$$
C_n^2 = |I|\cdot |J| \le |I\cap J|\cdot |L| = C_{n-1} C_{n+1}
$$
finishing the proof.
\eprf

We now consider the Motzkin numbers.  A {\em Motzkin path of length $n$} is a lattice path $P$ which satisfies
\ben
\item $P$ starts at $p_0=(0,0)$ and ends at $p_n=(n,0)$,
\item $P$ uses steps $U$, $D$, and horizontal $H=[1,0]$  and never goes below the $x$-axis.
\een
Let 
$$
\cM_n =\{P \mid \text{$P$ is a Motzkin path of length $n$}\}
$$
so that
$$
|\cM_n| = M_n,
$$
the $n$th {\em Motzkin number}. 
The set $\cM_n$ ordered by~\eqref{AP} is a distributive lattice as demonstrated in~\cite{FP:llp}.
Showing that the Motzkin sequence is log-convex is much like the proof  of the previous theorem, with $I$ and $J$ replaced by
$$
I=\{P\in\cM_{n+1} \mid \text{$P=HP'$ for some translated Motzkin  path $P'$ of length $n$}\},
$$
and
$$
J=\{P\in\cM_{n+1} \mid \text{$P=P'H$ for some  Motzkin  path $P'$ of length $n$}\}.
$$
So, we leave the details to the reader.  This result was proved algebraically by Aigner~\cite{aig:mn}.  A combinatorial proof for both Motzkin and the large Schr\"oder numbers was given by Chen, Wang, and Zheng~\cite{CWZ:cpl}.
\bth
The sequence $(M_n)_{n\ge0}$ of Motzkin numbers is log-convex.\hqed
\eth

Finally, we investigate the large Schr\"oder numbers.
A {\em Schr\"oder path of semilength $n$} is a lattice path $P$ satisfying
\ben
\item $P$ starts at $p_0=(0,0)$ and ends at $p_n=(n,0)$,
\item $P$ uses steps $U$, $D$, and twice horizontal $T=[2,0]$  and never goes below the $x$-axis.
\een
If we let
$$
\cS_n = \{P \mid \text{$P$ is a Schr\"oder path of semilength $n$}\}
$$
then
$$
|\cS_n| = S_n,
$$
the $n$th {\em large Schr\"oder number}.  As usual, we order $\cS_n$ using~\eqref{AP}.  However, this poset is not covered by the general theorem of~\cite{FP:llp}, although they remark that it can be shown that the poset is a lattice.  It is, in fact, distributive.
\ble
The poset $\cS_n$ is a distributive lattice.
\ele
\bprf
We first wish to show the existence of meets and joins.  Suppose $K,L\in \cS_n$.  We claim that there is $K\mt L\in\cS_n$ such that 
$$
A(K\mt L) = A(K)\cap A(L).
$$
To see this, consider any Schr\"oder path $L:p_0,p_1,\ldots,p_n$.  Since $L$ starts at $(0,0)$ and uses steps $U,H,T$ it must be that the coordinates of each $p_i$  have the same parity.
It follows that $K$ and $L$ can only intersect at points which are endpoints of steps in both paths and cannot cross in the middle of a step.  From this, the assertion that $A(K)\cap A(L)$ is the area under some Schr\"oder path follows.  Similarly, there is $K\jn L$ defined by
$$
A(K\jn L) = A(K)\cup A(L).
$$
It is now an easy matter to verify that $K\mt L$ and $K\jn L$ are indeed greatest lower bounds and least upper bounds.  And distributivity follows from the fact that intersection distributes over union.
\eprf

Using the ideals
$$
I=\{P\in\cS_{n+1} \mid \text{$P=TP'$ for some translated Schr\"oder  path $P'$ of semilength $n$}\},
$$
and
$$
J=\{P\in\cS_{n+1} \mid \text{$P=P'T$ for some  Schr\"oder  path $P'$ of semilength $n$}\},
$$
the next result follows in the way to which we have become accustomed.
\bth
The sequence $(S_n)_{n\ge0}$ of Schr\"oder  numbers is log-convex.\hqed
\eth

\section{Young's lattice}
\label{yl}

In this section we collect various log-concavity and  log-convexity results which can be proved using closed  intervals in Young's lattice, $Y$.  Since $Y$ is a distributive lattice, so is every closed interval and we need not check that part of the hypothesis of the Order Ideal Lemma.  We start by establishing some standard definitions and  notation, including
$$
[n]=\{1,2,\ldots,n\}.
$$

An {\em integer partition} is a weakly decreasing sequence $\la=(\la_1,\la_2,\ldots,\la_k)$ of nonnegative integers called {\em parts}.  Note that we allow zero as a part.
The {\em Young diagram} of $\la$ consists of left-justified rows of boxes or cells with $\la_i$ boxes in row $i$. 
We will use English notation for Young diagrams with the first row on top and make no distinction between a partition and its Young diagram.
For example, the shape of $\la=(5,3,1)$ is shown in the upper left in Figure~\ref{531}.
{\em Young's lattice}, $Y$, is the set of all partitions with partial order $\la\preceq\mu$ if the Young diagram of $\la$ is contained in that of $\mu$.

We first show how our proof of the log-convexity of the Catalan numbers can be viewed in this setting.  Consider the interval 
$$
\cC_n = [\emp,(n-1,n-2,\ldots,0)]
$$
 The southeast boundary of any $\la\in\cC_n$ can be considered as a lattice path of length $2n$ from the lower left corner to the upper right corner of  $(n-1,n-2,\ldots,0)$ using the edges of squares as unit steps north or east.  And the  lattice path must stay above the line $y=x$ if its initial point is taken as $(0,0)$.  Such paths are another of the standard ways of describing Dyck paths.  So the proof of Theorem~\ref{CnThm} could have been given in the language of partitions.

For our next result, if $\la=(\la_1,\la_2,\ldots,\la_k)$ is a partition and $n0$  is an integer then we define
\beq
\label{la+n}
\la+n = (\la_1+n,\la_2+n,\ldots,\la_k+n).
\eeq
\bth
\label{laThm}
For any partition $\la$, the sequence $(\ \#[\la,\la+n]\ )_{n\ge0}$ is log-concave.
\eth
\bprf
Let $I_n = [\la,\la+n]$ and consider the following as  subsets of $I_{n+1}$:
$$
I =[\la,\la+n]
$$
and
$$
J = [\la+1,\la+n+1].
$$
As in the previous proof, $I$ is an lower order ideal, $J$ is an upper order ideal, and $I=I_n\iso J$.
As expected, $I\cap J = [\la+1,\la+n]\iso I_{n-1}$.  So we are done by the Order Ideal Lemma.
\eprf

Note that even more generally, $\la$ could be replaced by a skew partition in the previous result and the proof would go through without change.
As an application, suppose that $\la=(1^k)$ where we us the notation $i^m$  to indicate that $i$ is repeated $m$ times
Then the $\mu\in[\la,\la+n]$ are precisely the $k$-element multisubsets of $[n+1]$.  These are counted by $\binom{n+k}{k}$ and so we immediately get  the following specialization of Theorem~\ref{laThm}.
\bco
\label{BinCo}
For any $k\ge0$, the sequence $\left(\ \binom{n+k}{k}\ \right)_{n\ge0}$ is log-concave. \hqed
\eco

\section{Order polynomials}
\label{op}

\bfi
\bce
\begin{tikzpicture}
\draw(0,0) circle(.2);
\draw(0,0) node{$2$};
\draw(-1,1) circle(.2);
\draw(-1,1) node{$1$};
\draw(1,1) circle(.2);
\draw(1,1) node{$3$};
\draw (-.15,.15)--(-.85,.85);
\draw (.15,.15)--(.85,.85);
\end{tikzpicture}
\ece 
\capt{A poset $P$ on $[3]$}
\label{P[3]}
\efi

The order polynomial of a labeled poset was introduced by Stanley in his thesis~\cite{sta:osp} and has since been shown to be a fundamental invariant.  In this section we prove that the sequence of values of the order polynomial of any labeled poset is log-concave.  This permits us to prove log-concavity of sequences formed by specializing the Schur function corresponding to any partition.  The enriched order polynomials of Stembridge~\cite{ste:epp} are also shown to give rise to log-concave sequences, again independent of the underlying poset.  
This gives rise to log-concave sequences of specializations of Schur $Q$-functions.
We begin, as usual, with the necessary definitions.  Note that some of the inequalities have been reversed from Stanley's original definitions, but this does not change the theory in any substantive way.

Let $(P,\peq)$ be a poset on $[p]$.  The reader should be sure to distinguish the use of $\preceq$ for the partial order on $P$ and $\leq$ for the total order on the integers.
A poset on $[3]$ is displayed in Figure~\ref{P[3]}.
A {\em $P$-partition with range $[n]$} is a map $f:P\ra[n]$ such that for all $x\prec y$ we have
\ben
\item $f(x)\le f(y)$ (that is, $f$ is order preserving), and
\item if $x>y$ then $f(x)<f(y)$.
\een
Define
$$
\cO_P(n) = \{ f \mid \text{$f$ is a  $P$-partition with range $n$}\}.
$$
Returning to the example poset in Figure~\ref{P[3]} we have 
$$
\cO_P(n) = \{f:P\ra[n] \mid f(2)<f(1) \text{ and } f(2)\le f(3)\}.
$$
The {\em order polynomial} of $P$ is
$$
\Om_P(n) = \# \cO_P(n).
$$
\bth[\cite{sta:osp}]
For any $P$ on $[p]$ we have $\Om_P(n)$ is a polynomial in $n$.\hqed
\eth

\bfi
\bce
\begin{tikzpicture}[scale=1]
\draw (-2,1.2)--(-3,1.8)  (2,3.2)--(.7,4) (6,5.2)--(4.7,6) (-2,5.2)--(-3,5.8) (2,7.2)--(.7,8)
(2,1.2)--(3,1.8) (6,3.2)--(7,3.8) (-2,3.2)--(-1,3.8) (2,5.2)--(3,5.8) (-2,7.2)--(-1,7.8) ;
\draw(0,0) circle(.2);
\node[scale = .5] at (0,0) {$2$};
\draw(-1,1) circle(.2);
\node[scale = .5] at (-1,1) {$1$};
\draw(1,1) circle(.2);
\node[scale = .5] at (1,1) {$3$};
\draw (-.15,.15)--(-.85,.85);
\draw (.15,.15)--(.85,.85);
\draw(-1.5,1) node{$2$};
\draw(1.5,1) node{$1$};
\draw(-.5,0) node{$1$};
\begin{scope}[shift={(-4,2)}]
\draw(0,0) circle(.2);
\node[scale = .5] at (0,0) {$2$};
\draw(-1,1) circle(.2);
\node[scale = .5] at (-1,1) {$1$};
\draw(1,1) circle(.2);
\node[scale = .5] at (1,1) {$3$};
\draw (-.15,.15)--(-.85,.85);
\draw (.15,.15)--(.85,.85);
\draw(-1.5,1) node{$3$};
\draw(1.5,1) node{$1$};
\draw(-.5,0) node{$1$};
\end{scope}
\begin{scope}[shift={(4,2)}]
\draw(0,0) circle(.2);
\node[scale = .5] at (0,0) {$2$};
\draw(-1,1) circle(.2);
\node[scale = .5] at (-1,1) {$1$};
\draw(1,1) circle(.2);
\node[scale = .5] at (1,1) {$3$};
\draw (-.15,.15)--(-.85,.85);
\draw (.15,.15)--(.85,.85);
\draw(-1.5,1) node{$2$};
\draw(1.5,1) node{$2$};
\draw(-.5,0) node{$1$};
\end{scope}
\begin{scope}[shift={(0,4)}]
\draw(0,0) circle(.2);
\node[scale = .5] at (0,0) {$2$};
\draw(-1,1) circle(.2);
\node[scale = .5] at (-1,1) {$1$};
\draw(1,1) circle(.2);
\node[scale = .5] at (1,1) {$3$};
\draw (-.15,.15)--(-.85,.85);
\draw (.15,.15)--(.85,.85);
\draw(-1.5,1) node{$3$};
\draw(1.5,1) node{$2$};
\draw(-.5,0) node{$1$};
\end{scope}
\begin{scope}[shift={(8,4)}]
\draw(0,0) circle(.2);
\node[scale = .5] at (0,0) {$2$};
\draw(-1,1) circle(.2);
\node[scale = .5] at (-1,1) {$1$};
\draw(1,1) circle(.2);
\node[scale = .5] at (1,1) {$3$};
\draw (-.15,.15)--(-.85,.85);
\draw (.15,.15)--(.85,.85);
\draw(-1.5,1) node{$2$};
\draw(1.5,1) node{$3$};
\draw(-.5,0) node{$1$};
\end{scope}
\begin{scope}[shift={(-4,6)}]
\draw(0,0) circle(.2);
\node[scale = .5] at (0,0) {$2$};
\draw(-1,1) circle(.2);
\node[scale = .5] at (-1,1) {$1$};
\draw(1,1) circle(.2);
\node[scale = .5] at (1,1) {$3$};
\draw (-.15,.15)--(-.85,.85);
\draw (.15,.15)--(.85,.85);
\draw(-1.5,1) node{$3$};
\draw(1.5,1) node{$2$};
\draw(-.5,0) node{$2$};
\end{scope}
\begin{scope}[shift={(4,6)}]
\draw(0,0) circle(.2);
\node[scale = .5] at (0,0) {$2$};
\draw(-1,1) circle(.2);
\node[scale = .5] at (-1,1) {$1$};
\draw(1,1) circle(.2);
\node[scale = .5] at (1,1) {$3$};
\draw (-.15,.15)--(-.85,.85);
\draw (.15,.15)--(.85,.85);
\draw(-1.5,1) node{$3$};
\draw(1.5,1) node{$3$};
\draw(-.5,0) node{$1$};
\end{scope}
\begin{scope}[shift={(0,8)}]
\draw(0,0) circle(.2);
\node[scale = .5] at (0,0) {$2$};
\draw(-1,1) circle(.2);
\node[scale = .5] at (-1,1) {$1$};
\draw(1,1) circle(.2);
\node[scale = .5] at (1,1) {$3$};
\draw (-.15,.15)--(-.85,.85);
\draw (.15,.15)--(.85,.85);
\draw(-1.5,1) node{$3$};
\draw(1.5,1) node{$3$};
\draw(-.5,0) node{$2$};
\end{scope}
\end{tikzpicture}
\ece
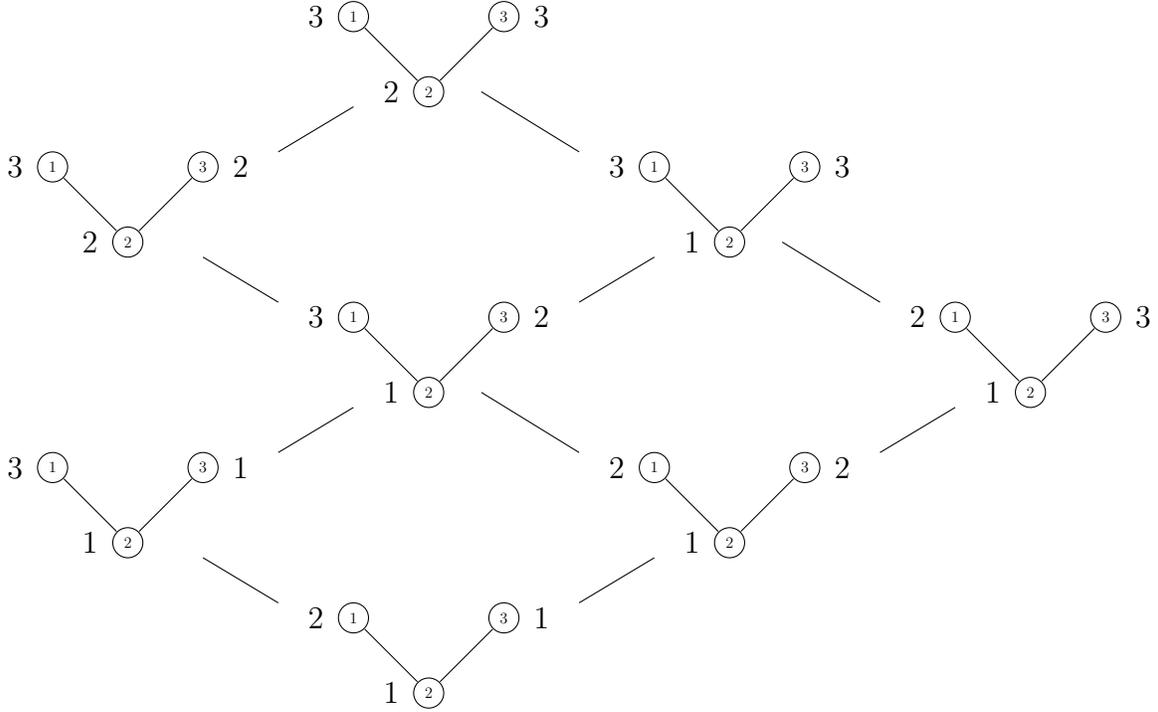
\capt{The poset $\cO_P(3)$ where $P$ is in Figure~\ref{P[3]}}
\label{cO_P:fig}
\efi

We now turn $\cO_p(n)$ into a poset by ordering $P$-partitions component-wise, that is,
$$
\text{$f\le g$  if and only if  $f(x)\le g(x)$ for all $x\in P$.}
$$
Even though we are using $\le$ for both the partial order on functions and the total order on integers, context  should distinguish them.
Continuing to use the poset in Figure~\ref{P[3]}, we have displayed the partial order on $\cO_P(3)$ in Figure~\ref{cO_P:fig}.  Note that the values of the functions $f$ are displayed outside the circles containing the elements of $P$ itself.
The following result was proved in the special case that $P$ is naturally labeled (that is, 
$x\prec y$ implies $x< y$) by Chan, Pak and Panova~\cite{CPP:epi}, and in full generality by Brenti~\cite[Theorem 7.6.5]{bre:ulc}
\bth
\label{Om:thm}
For any $P$ on $[p]$, the sequence $(\Om_P(n))_{n\ge0}$ is log-concave.
\eth
\bprf
We first need to prove that $\cO_P(n)$ is a distributive lattice.  Given $f,g\in\cO_P(n)$ we define a new function $f\mt g:P\ra[n]$ by taking component-wise minima
$$
(f\mt g)(x) = \min\{f(x),g(x)\}
$$
for all $x\in P$.  This will clearly be a greatest lower bound provided that $f\mt g$ is a $P$-partition.  Checking the two axioms are similar, so we only do the second.  So suppose that 
$x\peq y$ and $x>y$.  Then $f(x)<f(y)$ and $g(x)<g(y)$.  Without loss of generality we can assume that $f(y)\le g(y)$.  So,
$$
\min\{f(x),g(x)\} \le f(x) < f(y) =\min\{f(y),g(y)\}
$$
as desired.

Similarly, we define $f\jn g$ by
$$
(f\jn g)(x) = \max\{f(x),g(x)\}
$$
for all $x\in P$.  The fact that this is a least upper bound is much like what was done for the meet.  And the distributive law follows from the fact that maximum distributes over minimum.

To apply the Order Ideal Lemma, let
$$
L=\cO_P(n+1)
$$
so that $|L|=\Om_P(n+1)$.  Now define
$$
I = \{f\in L \mid f(x)\le n \text{ for all } x\in P\}.
$$
Clearly $I$ is a lower order ideal since we are applying an upper bound on $f$.  Furthermore, $I=\cO_P(n)$ so that $|I|=\Om_p(n)$.  Also consider
$$
J = \{f\in L \mid f(x)\ge 2 \text{ for all } x\in P\}.
$$
This is an upper order ideal and $J\iso\cO_P(n)$ where the isomorphism is obtained by subtracting one from each $f(x)$.  It follows that $|J|=\Om_P(n)$.  Similar considerations show that
$I\cap J\iso \cO_P(n-1)$ so that $|I\cap J|=\Om_P(n-1)$.  This completes the proof.
\eprf

\bfi
$$
\barr{cc}
\la=\ \ydiagram{5,3,1}
&
\qquad
\qquad
T =\ 
\begin{ytableau}
1&1&1&2&4\\
2&2&5\\
4
\end{ytableau}
\\[50pt]
\begin{tikzpicture}
\draw(-1,2) node{$P_\la=$}; 
\draw (0,4)--(4,0) (5,1)--(3,3) (4,0)--(6,2) (3,1)--(4,2) (2,2)--(3,3);
\filldraw[color=white](0,4) circle(.2);
\draw(0,4) circle(.2);
\draw(0,4) node{$9$};
\filldraw[color=white](1,3) circle(.2);
\draw(1,3) circle(.2);
\draw(1,3) node{$8$};
\filldraw[color=white](2,2) circle(.2);
\draw(2,2) circle(.2);
\draw(2,2) node{$7$};
\filldraw[color=white](3,1) circle(.2);
\draw(3,1) circle(.2);
\draw(3,1) node{$6$};
\filldraw[color=white](4,0) circle(.2);
\draw(4,0) circle(.2);
\draw(4,0) node{$5$};
\filldraw[color=white](5,1) circle(.2);
\draw(5,1) circle(.2);
\draw(5,1) node{$2$};
\filldraw[color=white](4,2) circle(.2);
\draw(4,2) circle(.2);
\draw(4,2) node{$3$};
\filldraw[color=white](3,3) circle(.2);
\draw(3,3) circle(.2);
\draw(3,3) node{$4$};
\filldraw[color=white](6,2) circle(.2);
\draw(6,2) circle(.2);
\draw(6,2) node{$1$};
\end{tikzpicture}
&
\begin{tikzpicture}
\draw (0,4)--(4,0) (5,1)--(3,3) (4,0)--(6,2) (3,1)--(4,2) (2,2)--(3,3);
\filldraw[color=white](0,4) circle(.2);
\draw(0,4) circle(.2);
\draw(0,4) node{$9$};
\draw(-.5,4) node{$4$};
\filldraw[color=white](1,3) circle(.2);
\draw(1,3) circle(.2);
\draw(1,3) node{$8$};
\draw(.5,3) node{$2$};
\filldraw[color=white](2,2) circle(.2);
\draw(2,2) circle(.2);
\draw(2,2) node{$7$};
\draw(1.5,2) node{$1$};
\filldraw[color=white](3,1) circle(.2);
\draw(3,1) circle(.2);
\draw(3,1) node{$6$};
\draw(2.5,1) node{$1$};
\filldraw[color=white](4,0) circle(.2);
\draw(4,0) circle(.2);
\draw(4,0) node{$5$};
\draw(3.5,0) node{$1$};
\filldraw[color=white](5,1) circle(.2);
\draw(5,1) circle(.2);
\draw(5,1) node{$2$};
\draw(4.5,1) node{$2$};
\filldraw[color=white](4,2) circle(.2);
\draw(4,2) circle(.2);
\draw(4,2) node{$3$};
\draw(3.5,2) node{$2$};
\filldraw[color=white](3,3) circle(.2);
\draw(3,3) circle(.2);
\draw(3,3) node{$4$};
\draw(2.5,3) node{$5$};
\filldraw[color=white](6,2) circle(.2);
\draw(6,2) circle(.2);
\draw(6,2) node{$1$};
\draw(5.5,2) node{$4$};
\end{tikzpicture}
\earr
$$
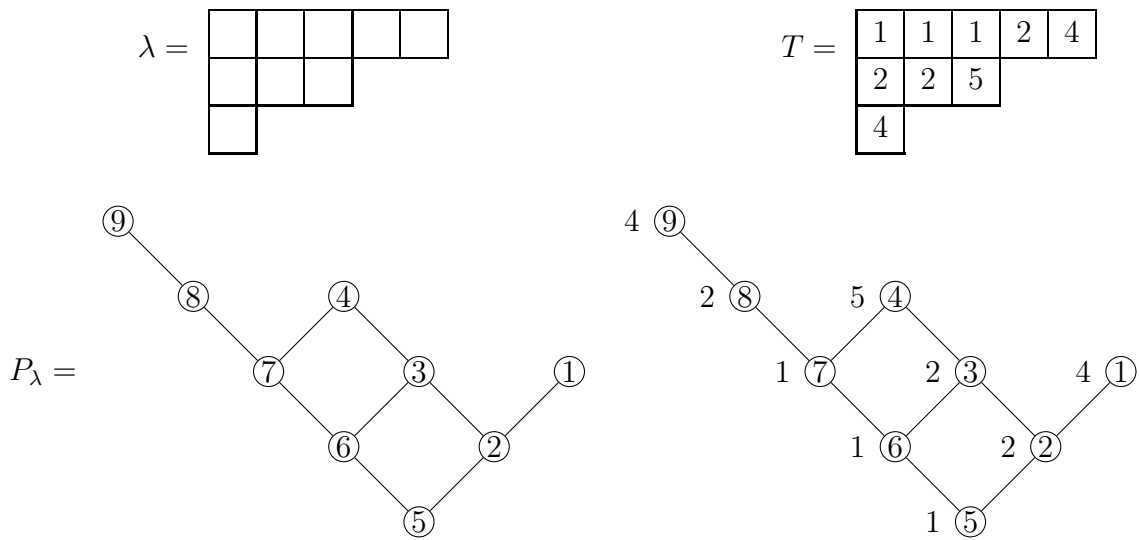
\capt{The shape $\la=(5,3,1)$, a semistandard Young tableau, $T$ of that shape, as well as the corresponding poset $P_\la$ and $P$-partition.}
\label{531}
\efi

We now use the well-known connection between order polynomials and Schur functions to derive an interesting special case of the previous theorem.  If $\la$ is an integer partition then a {\em semistandard Young tableau (SSYT) of shape $\la$} is a filling of the boxes of $\la$ with positive integers such that rows weakly increase left-to-right  and columns strictly increase top-to-bottom.  
The partition $\la=(5,3,1)$ and a semistandard Young tableau $T$ of that shape are displayed in the first row of Figure~\ref{531}.
We let $(i,j)$ be the cell of $\la$ in row $i$ and column $j$ where rows and columns are indexed as in a matrix. 
Given an SSYT of shape $\la$ we denote by $T_{i,j}$ the element of $T$ in box $(i,j)$.
In the tableau of Figure~\ref{531} we have $T_{2,3}=5$.
 Consider
$$
\SSYT_\la = \{T \mid \text{$T$ is an SSYT of shape $\la$}\}.
$$
Let  $\bx=\{x_1,x_2,\ldots\}$ be a set of variables indexed by the positive integers.  The {\em Schur function} corresponding to $\la$ is the generating function
$$
s_\la(\bx) = \sum_{T\in \SSYT_\la} \prod_{(i,j)\in\la} x_{T_{i,j}}.
$$
The Schur functions are symmetric and form an important basis for the algebra of symmetric functions.  For more information about them, see the texts of Sagan~\cite{sag:sg} or Stanley~\cite{sta:ec2}.

To make the connection with $P$-partitions, we first turn $\la=(\la_1,\la_2,\ldots,\la_k)$ into a poset  component-wise, that is
\beq
\label{Pla}
\text{$(i,j) \peq (i',j')$  if and only if $i\le i'$ and $j\le j'$.} 
\eeq
We now make this a poset $P_\la$ on the interval $[|\la|]$ where $|\la|=\sum_l \la_l$ by labeling the last row of $\la$ with $1,2,\ldots,\la_k$ left-to-right (viewing $\la$ as its original Young diagram).
Then labeling the penultimate row left-to-right with $\la_k+1,\la_k+2,\ldots,\la_k+\la_{k-1}$, and so forth.
This labeling is displayed in Figure~\ref{531} at the bottom left.  It is easy to see that a $P_\la$-partition is the same as an SSYT of shape $\la$.  The partition for the SSYT $T$ in Figure~\ref{531} is displayed directly below the tableau.  It should now be clear that we have 
$$
s_\la(1^n) = \Om_{P_\la}(n)
$$
where $1^n$ indicates the specialization
$$
\text{$x_i=1$ for $i\le n$ and $x_i=0$ for $i>n$.}
$$
As an immediate consequence of Theorem~\ref{Om:thm}, we have the following result which is also implicit in Brenti's work~\cite[Theorems 5.2.3 and 7.6.5]{bre:ulc}.
\bco
\label{sla}
For any partition $\la$, the sequence $(s_\la(1^n))_{n\ge0}$ is log-concave.\hqed
\eco

\bfi
\bce
\begin{tikzpicture}[scale=.9]
\draw (-2,1.2)--(-3,1.8) (2,3.2)--(.7,4) (-6,5.2)--(-7,5.8) (-6,7.2)--(-7,7.8) (-2,7.2)--(-3,7.8) (-1.8,9.2)--(-3,9.8) (2.2,9.2)--(1,9.8) (2,11.2)--(1,11.8)
(2,1.2)--(3,1.8) (-2,3.2)--(-1,3.8) (-2,5.2)--(-1,5.8) (2,5.2)--(3,5.8) (2,7.2)--(3,7.8) (-6,9.2)--(-5,9.8) (-2.2,9.2)--(-1,9.8) (1.8,9.2)--(3,9.8) (-2,11.2)--(-1,11.8)
(-4,3)--(-4,3.6) (4,3)--(4,3.6) (-4,5)--(-4,5.6) (0,5)--(0,5.6) (4,5)--(4,5.6) (-8,7)--(-8,7.6) (-4,7)--(-4,7.6) (0,7)--(0,7.6) (4,7)--(4,7.6) (-4,9)--(-4,9.6) (4,9)--(4,9.6) (0,11)--(0,11.6) (0,13)--(0,13.6)
(-6.5,7.3)--(-.7,8);
\draw(0,0) circle(.2);
\node[scale = .5] at (0,0) {$2$};
\draw(-1,1) circle(.2);
\node[scale = .5] at (-1,1) {$1$};
\draw(1,1) circle(.2);
\node[scale = .5] at (1,1) {$3$};
\draw (-.15,.15)--(-.85,.85);
\draw (.15,.15)--(.85,.85);
\draw(-1.5,1) node{$\ol{1}$};
\draw(1.5,1) node{$1$};
\draw(-.5,0) node{$\ol{1}$};
\begin{scope}[shift={(-4,2)}]
\draw(0,0) circle(.2);
\node[scale = .5] at (0,0) {$2$};
\draw(-1,1) circle(.2);
\node[scale = .5] at (-1,1) {$1$};
\draw(1,1) circle(.2);
\node[scale = .5] at (1,1) {$3$};
\draw (-.15,.15)--(-.85,.85);
\draw (.15,.15)--(.85,.85);
\draw(-1.5,1) node{$1$};
\draw(1.5,1) node{$1$};
\draw(-.5,0) node{$\ol{1}$};
\end{scope}
\begin{scope}[shift={(4,2)}]
\draw(0,0) circle(.2);
\node[scale = .5] at (0,0) {$2$};
\draw(-1,1) circle(.2);
\node[scale = .5] at (-1,1) {$1$};
\draw(1,1) circle(.2);
\node[scale = .5] at (1,1) {$3$};
\draw (-.15,.15)--(-.85,.85);
\draw (.15,.15)--(.85,.85);
\draw(-1.5,1) node{$\ol{1}$};
\draw(1.5,1) node{$\ol{2}$};
\draw(-.5,0) node{$\ol{1}$};
\end{scope}
\begin{scope}[shift={(-4,4)}]
\draw(0,0) circle(.2);
\node[scale = .5] at (0,0) {$2$};
\draw(-1,1) circle(.2);
\node[scale = .5] at (-1,1) {$1$};
\draw(1,1) circle(.2);
\node[scale = .5] at (1,1) {$3$};
\draw (-.15,.15)--(-.85,.85);
\draw (.15,.15)--(.85,.85);
\draw(-1.5,1) node{$\ol{2}$};
\draw(1.5,1) node{$1$};
\draw(-.5,0) node{$\ol{1}$};
\end{scope}
\begin{scope}[shift={(0,4)}]
\draw(0,0) circle(.2);
\node[scale = .5] at (0,0) {$2$};
\draw(-1,1) circle(.2);
\node[scale = .5] at (-1,1) {$1$};
\draw(1,1) circle(.2);
\node[scale = .5] at (1,1) {$3$};
\draw (-.15,.15)--(-.85,.85);
\draw (.15,.15)--(.85,.85);
\draw(-1.5,1) node{$1$};
\draw(1.5,1) node{$\ol{2}$};
\draw(-.5,0) node{$\ol{1}$};
\end{scope}
\begin{scope}[shift={(4,4)}]
\draw(0,0) circle(.2);
\node[scale = .5] at (0,0) {$2$};
\draw(-1,1) circle(.2);
\node[scale = .5] at (-1,1) {$1$};
\draw(1,1) circle(.2);
\node[scale = .5] at (1,1) {$3$};
\draw (-.15,.15)--(-.85,.85);
\draw (.15,.15)--(.85,.85);
\draw(-1.5,1) node{$\ol{1}$};
\draw(1.5,1) node{$2$};
\draw(-.5,0) node{$\ol{1}$};
\end{scope}
\begin{scope}[shift={(-8,6)}]
\draw(0,0) circle(.2);
\node[scale = .5] at (0,0) {$2$};
\draw(-1,1) circle(.2);
\node[scale = .5] at (-1,1) {$1$};
\draw(1,1) circle(.2);
\node[scale = .5] at (1,1) {$3$};
\draw (-.15,.15)--(-.85,.85);
\draw (.15,.15)--(.85,.85);
\draw(-1.5,1) node{$2$};
\draw(1.5,1) node{$1$};
\draw(-.5,0) node{$\ol{1}$};
\end{scope}
\begin{scope}[shift={(-4,6)}]
\draw(0,0) circle(.2);
\node[scale = .5] at (0,0) {$2$};
\draw(-1,1) circle(.2);
\node[scale = .5] at (-1,1) {$1$};
\draw(1,1) circle(.2);
\node[scale = .5] at (1,1) {$3$};
\draw (-.15,.15)--(-.85,.85);
\draw (.15,.15)--(.85,.85);
\draw(-1.5,1) node{$\ol{2}$};
\draw(1.5,1) node{$1$};
\draw(-.5,0) node{$1$};
\end{scope}
\begin{scope}[shift={(0,6)}]
\draw(0,0) circle(.2);
\node[scale = .5] at (0,0) {$2$};
\draw(-1,1) circle(.2);
\node[scale = .5] at (-1,1) {$1$};
\draw(1,1) circle(.2);
\node[scale = .5] at (1,1) {$3$};
\draw (-.15,.15)--(-.85,.85);
\draw (.15,.15)--(.85,.85);
\draw(-1.5,1) node{$\ol{2}$};
\draw(1.5,1) node{$\ol{2}$};
\draw(-.5,0) node{$\ol{1}$};
\end{scope}
\begin{scope}[shift={(4,6)}]
\draw(0,0) circle(.2);
\node[scale = .5] at (0,0) {$2$};
\draw(-1,1) circle(.2);
\node[scale = .5] at (-1,1) {$1$};
\draw(1,1) circle(.2);
\node[scale = .5] at (1,1) {$3$};
\draw (-.15,.15)--(-.85,.85);
\draw (.15,.15)--(.85,.85);
\draw(-1.5,1) node{$1$};
\draw(1.5,1) node{$2$};
\draw(-.5,0) node{$\ol{1}$};
\end{scope}
\begin{scope}[shift={(-8,8)}]
\draw(0,0) circle(.2);
\node[scale = .5] at (0,0) {$2$};
\draw(-1,1) circle(.2);
\node[scale = .5] at (-1,1) {$1$};
\draw(1,1) circle(.2);
\node[scale = .5] at (1,1) {$3$};
\draw (-.15,.15)--(-.85,.85);
\draw (.15,.15)--(.85,.85);
\draw(-1.5,1) node{$2$};
\draw(1.5,1) node{$1$};
\draw(-.5,0) node{$1$};
\end{scope}
\begin{scope}[shift={(-4,8)}]
\draw(0,0) circle(.2);
\node[scale = .5] at (0,0) {$2$};
\draw(-1,1) circle(.2);
\node[scale = .5] at (-1,1) {$1$};
\draw(1,1) circle(.2);
\node[scale = .5] at (1,1) {$3$};
\draw (-.15,.15)--(-.85,.85);
\draw (.15,.15)--(.85,.85);
\draw(-1.5,1) node{$\ol{2}$};
\draw(1.5,1) node{$\ol{2}$};
\draw(-.5,0) node{$1$};
\end{scope}
\begin{scope}[shift={(0,8)}]
\draw(0,0) circle(.2);
\node[scale = .5] at (0,0) {$2$};
\draw(-1,1) circle(.2);
\node[scale = .5] at (-1,1) {$1$};
\draw(1,1) circle(.2);
\node[scale = .5] at (1,1) {$3$};
\draw (-.15,.15)--(-.85,.85);
\draw (.15,.15)--(.85,.85);
\draw(-1.5,1) node{$2$};
\draw(1.5,1) node{$\ol{2}$};
\draw(-.5,0) node{$\ol{1}$};
\end{scope}
\begin{scope}[shift={(4,8)}]
\draw(0,0) circle(.2);
\node[scale = .5] at (0,0) {$2$};
\draw(-1,1) circle(.2);
\node[scale = .5] at (-1,1) {$1$};
\draw(1,1) circle(.2);
\node[scale = .5] at (1,1) {$3$};
\draw (-.15,.15)--(-.85,.85);
\draw (.15,.15)--(.85,.85);
\draw(-1.5,1) node{$\ol{2}$};
\draw(1.5,1) node{$2$};
\draw(-.5,0) node{$\ol{1}$};
\end{scope}
\begin{scope}[shift={(-4,10)}]
\draw(0,0) circle(.2);
\node[scale = .5] at (0,0) {$2$};
\draw(-1,1) circle(.2);
\node[scale = .5] at (-1,1) {$1$};
\draw(1,1) circle(.2);
\node[scale = .5] at (1,1) {$3$};
\draw (-.15,.15)--(-.85,.85);
\draw (.15,.15)--(.85,.85);
\draw(-1.5,1) node{$2$};
\draw(1.5,1) node{$\ol{2}$};
\draw(-.5,0) node{$1$};
\end{scope}
\begin{scope}[shift={(0,10)}]
\draw(0,0) circle(.2);
\node[scale = .5] at (0,0) {$2$};
\draw(-1,1) circle(.2);
\node[scale = .5] at (-1,1) {$1$};
\draw(1,1) circle(.2);
\node[scale = .5] at (1,1) {$3$};
\draw (-.15,.15)--(-.85,.85);
\draw (.15,.15)--(.85,.85);
\draw(-1.5,1) node{$\ol{2}$};
\draw(1.5,1) node{$2$};
\draw(-.5,0) node{$1$};
\end{scope}
\begin{scope}[shift={(4,10)}]
\draw(0,0) circle(.2);
\node[scale = .5] at (0,0) {$2$};
\draw(-1,1) circle(.2);
\node[scale = .5] at (-1,1) {$1$};
\draw(1,1) circle(.2);
\node[scale = .5] at (1,1) {$3$};
\draw (-.15,.15)--(-.85,.85);
\draw (.15,.15)--(.85,.85);
\draw(-1.5,1) node{$2$};
\draw(1.5,1) node{$2$};
\draw(-.5,0) node{$\ol{1}$};
\end{scope}
\begin{scope}[shift={(0,12)}]
\draw(0,0) circle(.2);
\node[scale = .5] at (0,0) {$2$};
\draw(-1,1) circle(.2);
\node[scale = .5] at (-1,1) {$1$};
\draw(1,1) circle(.2);
\node[scale = .5] at (1,1) {$3$};
\draw (-.15,.15)--(-.85,.85);
\draw (.15,.15)--(.85,.85);
\draw(-1.5,1) node{$2$};
\draw(1.5,1) node{$2$};
\draw(-.5,0) node{$1$};
\end{scope}
\begin{scope}[shift={(0,14)}]
\draw(0,0) circle(.2);
\node[scale = .5] at (0,0) {$2$};
\draw(-1,1) circle(.2);
\node[scale = .5] at (-1,1) {$1$};
\draw(1,1) circle(.2);
\node[scale = .5] at (1,1) {$3$};
\draw (-.15,.15)--(-.85,.85);
\draw (.15,.15)--(.85,.85);
\draw(-1.5,1) node{$2$};
\draw(1.5,1) node{$2$};
\draw(-.5,0) node{$\ol{1}$};
\end{scope}
\end{tikzpicture}
\ece
\capt{The poset $\cO_P^e(2)$ where $P$ is in Figure~\ref{P[3]}}
\label{O_P^e}
\efi

We can also apply the Order Ideal Lemma to the enriched $P$-partitions of Stembridge~\cite{ste:epp}.  We  put an unusual total order on the nonzero integers, where we denote $-n$ by $\nb$,
\beq
\label{tles}
\ol{1} \tles 1 \tles \ol{2} \tles 2 \tles \ol{3} \tles 3 \tles \ldots\ 
\eeq
and also let
$$
\spn{n} = \{\ol{1}, 1, \ol{2}, 2, \ldots, \nb, n\}.
$$
An {\em enriched $P$-partition with range $\spn{n}$} is a map 
$f:P\ra\spn{n}$ such that for all $x\prec y$ we have
\ben
\item[(E1)] $f(x)\tleq f(y)$,
\item[(E2)] $f(x)=f(y)>0$ implies $x<y$, and
\item[(E3)] $f(x)=f(y)<0$ implies $x>y$.
\een
We now let
$$
\cO_P^e(n) = \{f \mid \text{$f$ is an enriched  $P$-partition with range $\spn{n}$}\}.
$$
with corresponding {\em enriched order polynomial}
$$
\Om_P^e(n) = \#\cO_P^e(n).
$$
As with the ordinary order polynomial, the enriched one is well named.
\bth[\cite{ste:epp}]
\label{eP:thm}
For any $P$ on $[p]$ we have $\Om_P^e(n)$ is a polynomial in $n$.\hqed
\eth

We now turn $\cO_P^e(n)$ into a poset in exactly the same way as the ordinary case:
$$
\text{$f\le g$  if and only if  $f(x)\tleq g(x)$ for all $x\in P$.}
$$
Figure~\ref{O_P^e} shows the partial order $\cO_P^e(2)$ using our canonical poset $P$ in Figure~\ref{P[3]}.
The proof of the next result is much like that of Theorem~\ref{Om:thm} and so is omitted.
\bth
\label{Ome:thm}
For any $P$ on $[p]$, the sequence $(\Om_P^e(n))_{n\ge0}$ is log-concave.\hqed
\eth

\bfi
\ydiagram{6,1+5,2+3,3+2}
\hspace{50pt}
\begin{ytableau}
\ol{1}  &   1   &   1   &\ol{2}&\ol{3}&   3   &   3\\
\none   &   2   &   2   &   2  &\ol{3}&\ol{4} &   4\\
\none   & \none &\ol{3} &   3  &    3\\
\none   & \none & \none &   4  &    4
\end{ytableau}
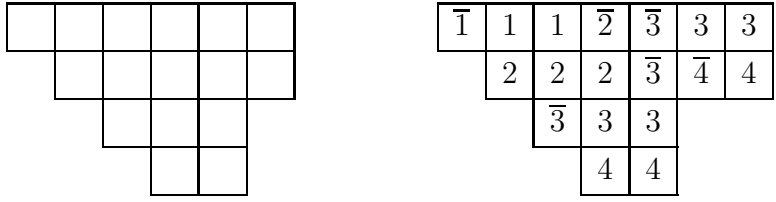
\capt{The shifted Young diagram of $\la=(6,5,3,2)$ and a shifted semistandard Young tableau of that shape}
\label{6532}
\efi

The appropriate tableaux for the enriched setting arise from shifted shapes.  An integer partition
$\la=(\la_1,\la_2,\ldots,\la_k)$ is {\em strict} if $\la_1>\la_2>\ldots>\la_k$.
A strict partition has an associated {\em shifted Young diagram} obtained from the ordinary Young diagram by shifting the $i$th row $i-1$ boxes to the right for $i\in[k]$.
The shifted shape of $\la=(6,5,3,2)$ is shown on the left in Figure~\ref{6532}.
A {\em semistandard Young tableau of shifted shape $\la$} is a filling $T$ of the cells of $\la$ with nonzero integers such that the following are satisfied.
\ben
\item[(T1)] The rows and columns of $T$ are weakly increasing with respect to the total order~\ref{tles}.
\item[(T2)] For each $m>0$ there is at most one $\ol{m}$ in each row and at most one $m$ in each column.
\een
We let
$$
\SSYTS_\la = \{T \mid \text{$T$ is a semistandard Young tableau of shifted shape $\la$}\}
$$
with generating function
$$
Q_\la(\bx) =\sum_{\SSYTS_\la} \prod_{(i,j)\in\la} x_{|T_{i,j}|}
$$
which is called a {\em Schur $Q$-function}.  These functions play a role in the projective representation theory of the symmetric group analogous to $s_\la(\bx)$ for ordinary representations.  

The partial order on shifted shape is just the restriction of~\eqref{Pla} to those cells in the shifted Young diagram.  And the labeling to get a corresponding partition $P_\la^e$ on $[|\la|]$ is exactly the same as in the unshifted case, starting with the bottom row and working up.  Now axiom (E1) for enriched $P$-partitions implies condition (T1) for semistandard shifted Young tableaux.  And axioms (E2) and (E3) translate into condition (T2).
Thus
$$
Q_\la(1^n) = \Om_{P_\la^e}(n)
$$
and the next result is a special case of Theorem~\ref{eP:thm}.
\bco
For any strict partition $\la$, the sequence $(Q_\la(1^n))_{n\ge0}$ is log-concave. \hqed
\eco

\section{Generalized Lucas sequences}
\label{gls}

A sequence $(l_n)_{n\ge0}$ of real numbers is a {\em generalized Lucas sequence} if it satisfies the recursion
\beq
\label{l:rr}
l_n=l_{n-1}+l_{n-2}
\eeq
for $n\ge2$.  These sequences were originally studied by Lucas~\cite{luc:tfn1,luc:tfn2,luc:tfn3}.  Both the sequences themselves and their $q$-analogues have many wonderful combinatorial properties, see~\cite{BP:caf,BCMS:cil,SS:cib,ST:la}.  Of course, the two most famous examples of such sequences are the {\em Fibonacci numbers}, $(F_n)_{n\ge0}$, and {\em (ordinary) Lucas numbers}, $(L_n)_{n\ge0}$, with initial conditions $F_0=F_1=1$ and $L_0=2$, $L_1=1$, respectively.

In this section we will study {\em positive Lucas sequences} which are generalized Lucas sequences with $l_0,l_1>0$.
In order to state our result precisely, call a sequence $(a_n)_{n\ge0}$ {\em log-concave at index $n$} if
$$
a_n^2 \ge a_{n-1} a_{n+1}.
$$
Note that this definition says nothing about indices other than $n$.  Similarly define {\em log-convexity at index $n$}. 
We will show that any positive Lucas sequence, suitably reindexed, alternates between being log-concave at odd indices and log-convex at even ones.

It will be convenient in our approach to restrict the initial values even further.  But we wish to first show that this restriction will be, in some sense, without loss of generality.  To do this, we extend a generalized Lucas sequence to negative indices by insisting that the recurrence relation~\eqref{l:rr} continue to hold for $n<0$ to give an {\em extended Lucas sequence} $(l_n)_{n\in\bbZ}$.  Call two extended Lucas sequences $(l_n)_{n\in\bbZ}$ and $(l'_n)_{n\in\bbZ}$ {\em shift equivalent} if there is $k\in\bbZ$ such that 
$$
l_n=l'_{n+k}
$$
for all $n\in\bbZ$.
\begin{prop}
Suppose that  $(l_n)_{n\ge0}$ is a positive Lucas sequence.  Then its extension is shift equivalent to an extended Lucas sequence $(l'_n)_{n\in\bbZ}$ such that
$$
0< 2l'_0 \le l'_1.
$$
\end{prop}
\begin{proof}
Consider the reverse subsequence $l_1,l_0,l_{-1},\ldots$ of the given Lucas sequence.  Suppose first that this sequence contains a weak ascent $l_j\le l_{j-1}$ where $l_j,l_{j-1}> 0$.  But then
$$
l_{j+1}=l_j+ l_{j-1}\ge 2 l_j >0,
$$
as desired. 

Now suppose this sequence contains a weak ascent with $l_j = 0$,  Then $l_{j-1}>0$ since otherwise all entries in the original Lucas sequence are nonpositive.
But now
$$
l_{j+2} = l_{j+1} + l_j = l_{j+1},
$$
where
$$
l_{j+1} = l_j + l_{j-1} = l_{j-1} >0.
$$
So
$$
l_{j+3} = l_{j+2} + l_{j+1} = 2 l_{j+2} >0,
$$
again giving the correct conclusion.

If there  are no such weak ascents, then the sequence is strictly decreasing and so must eventually become negative.  Let $m$ be the index of minimum absolute value such that $l_m<0$.  It follows that $l_{m+1}\ge0 $ and
$$
0\le l_{m+2} = l_{m+1}+ l_m\le l_{m+1}.
$$
But this is a weak ascent which is a contradiction.
\end{proof}

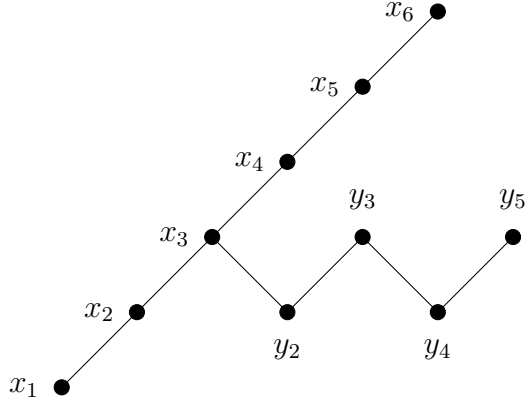
\begin{figure}
    \centering
\begin{tikzpicture}
\filldraw(1,1) circle(.1);
\draw(.5,1) node{$x_1$};
\filldraw(2,2) circle(.1);
\draw(1.5,2) node{$x_2$};
\filldraw(3,3) circle(.1);
\draw(2.5,3) node{$x_3$};
\filldraw(4,4) circle(.1);
\draw(3.5,4) node{$x_4$};
\filldraw(5,5) circle(.1);
\draw(4.5,5) node{$x_5$};
\filldraw(6,6) circle(.1);
\draw(5.5,6) node{$x_6$};
\filldraw(4,2) circle(.1);
\draw(4,1.5) node{$y_2$};
\filldraw(5,3) circle(.1);
\draw(5,3.5) node{$y_3$};
\filldraw(6,2) circle(.1);
\draw(6,1.5) node{$y_4$};
\filldraw(7,3) circle(.1);
\draw(7,3.5) node{$y_5$};
\draw (1,1)--(6,6) (3,3)--(4,2)--(5,3)--(6,2)--(7,3);
\end{tikzpicture}
    \caption{The poset $L_5(3,7)$}
    \label{L_5(3,7)}
\end{figure}

We will now introduce the posets whose lattices of order ideals will permit us to study the behaviour of positive Lucas sequences $(l_n)_{n\ge0}$.  Say that such a sequence is  {\em well-indexed} if
$$
0<2l_0\le l_1.
$$
Note that,  by the previous proposition, every positive Lucas sequence is shift equivalent to a well-indexed one.  To simplify notation, we will relabel
\begin{equation}
\label{rs}
 r:= l_0 \text{ and } s:= l_1   
 \end{equation}
Define a poset $L_n(r,s)$ to have elements $x_1,\ldots,x_{s-1}$ and $y_2,\ldots, y_n$ 
and order relation $\preceq$
subject to the covers
\begin{enumerate}
    \item $x_1\prec x_2 \prec \ldots \prec x_{s-1}$,
    \item $y_2 \prec y_3 \succ y_4\prec y_5 \succ \ldots$, and
    \item $y_2 \prec x_r$.
\end{enumerate}
So the $x_i$ form a chain $C_{s-1}$ and the $y_j$ form what we will call an {\em alternating poset} $A_{n-1}$.
For example, Figure~\ref{L_5(3,7)} shows the poset $L_5(3,7)$.

Let $J(P)$ denote the set of lower order ideals of a finite poset $P$.  It  is a fundamental result that $J(P)$ is a distributive lattice for any $P$.  This construction will permit us to prove the following theorem.  Note that since the Fibonacci sequence satisfies
$$
F_n^2 = F_{n-1} F_{n+1} + (-1)^{n-1},
$$
one sees immediately that the result is true in this case.
\begin{thm}
A well-indexed Lucas sequence $(l_n)_{n\ge0}$ is log-concave at odd indices and log-convex at even ones.
\end{thm}
\begin{proof}
We continue to use the notation in~\eqref{rs}.  We first verify log-concavity at index $n=1$ since our approach with ideals will only start to work when $n\ge2$.  Since $r\le s/2$ we have 
$$
l_0 l_2 = r (r +s) \le s/2(3s/2) <s^2 = l_1^2.
$$

For ease of notation, let $L_n=L_n(r,s)$.  We will show by induction that 
\begin{equation}
\label{Jnln}
    \#J(L_n) = l_n
\end{equation}
for $n\ge1$.  When $n=1$ we have that $L_1$ is a chain with 
$s-1$ elements so that $\#J(L_1)= s=l_1$.  Poset $L_2$ is obtained from $L_1$ by placing $y_2$ under the $r$th smallest element of the chain.  So, counting the ideals with $0,1,2,\ldots$ elements we get
$$
\#J(L_2) = 1 + \overbrace{2+\cdots+2}^{r-1}+\overbrace{1+\cdots+1}^{s-r+1}=r+s = l_2.
$$

The induction step has two cases depending on whether $n$ is odd or even.  But they are similar so we will only do the former.  In this case we have $y_{n-1}\prec y_n$.  Every ideal $I\in J(L_n)$ either contains $y_n$ or not.  It follows that if $y_n\not\in I$ then $I\in J(L_{n-1})$.
On the other hand, if $y_n\in I$ then this forces $y_{n-1}\in I$ and 
$I-\{y_{n-1},y_n\}\in J(L_{n-2})$.
From these observations and  induction
$$
\# J(L_n)=\# J(L_{n-1})+\# J(L_{n-2}) = l_{n-1}+l_{n-2} = l_n
$$
which completes the proof of the claim.

We now construct the ideals $\cI_n,\cJ_n\subseteq J(L_{n+1})$ needed to prove the theorem.  Given a set of constraints $S$ on ideals $I$ we will use the notation
$$
J_{n+1}(S) = \{I \in J(L_{n+1}) \mid \text{$I$ satisfies $S$}\}.
$$
Define
$$
\cI_n=
\case{J_{n+1}(y_{n+1}\in I)}{if $n$ is odd,}{J_{n+1}(y_{n+1}\not\in I)}{if $n$ is even.}
$$
Then
\begin{equation}
  \label{Inln}  
  \#\cI_n = l_n.
\end{equation}
In fact we have an isomorphism $\cI_n\iso J(L_n)$ given by the identity map on individual elements when $n$ is even,  and by removal of $y_{n+1}$ when $n$ is odd.  So the claim follows from~\eqref{Jnln}.
Note that the ideals in $\cI_n$ form either an upper or lower order ideal in $J(L_{n+1})$ depending on whether $n$ is odd or even, respectively.

Now let
$$
\cJ_n = J_{n+1}(x_r\not\in I) \hs{3pt} \cup \hs{3pt} J_{n+1}(x_{s-r}\not\in I  \text{ and } y_{n+1}\not\in I).
$$
Note that $\cJ_n$ is a lower order ideal in $J(L_{n+1})$ regardless of the parity of $n$.  We wish to show
\begin{equation}
\label{cJnln}
    \#\cJ_n = l_n.
\end{equation}
In view of~\eqref{Inln}, it suffices to show that the set differences $\cI_n\setm\cJ_n$ and $\cJ_n\setm\cI_n$ have the same cardinality.  Again, we merely provide details when $n$ is even.  Now
$$
\cI_n\setm\cJ_n = J_{n+1}(\text{$x_r\in I$, and $x_{s-r}\in I$, and $y_{n+1}\not\in I$}).
$$
But $r\le s-r$ so that $x_r\le x_{s-r}$ in $L_{n+1}$.  This makes the condition $x_r\in I$ redundant and
$$
\cI_n\setm\cJ_n = J_{n+1}(\text{$x_{s-r}\in I$, and $y_{n+1}\not\in I$}).
$$
Note that  $x_{s-r}$  is greater than $y_2$ and all the $x_i$ below it.  So those elements are forced to be in any ideal $I$ we are considering.  But then $I$ must be constructed by adding to these elements some ideal of the
chain $x_{s-r+1}, \ldots, x_{s-1}$  and, because $n$ is even, an ideal of the subposet of $L_{n+1}$ induced on $y_3,\ldots,y_n$.  It follows that we have a product poset
\begin{equation}
    \label{I-J}
    \cI_n\setm\cJ_n \iso J(C_{r-1}) \times J(A_{n-2}).
\end{equation}

Now directly from the definitions we have
$$
\cJ_n\setm\cI_n = J_{n+1}(\text{$x_r\not\in I$, and $y_{n+1}\in I$}).
$$
Since $x_r\not\in I$, any ideal in the difference breaks into two pieces.  One is an ideal in the
chain $x_1,\ldots,x_{r-1}$ and the other an ideal of the alternating poset on $y_2,\ldots,y_{n+1}$ which contains $y_{n+1}$.  Since by parity, $y_{n+1}\succ y_n$, we must also have $y_n\in I$ and so this part of $I$ is determined as an ideal in the alternating poset on $y_2,\ldots,y_{n-1}$.  Comparing this with~\eqref{I-J}
shows that $\#\cI_n\setm\cJ_n=\#\cJ_n\setm\cI_n$.

There remains to calculate $\#(\cI_n\cap \cJ_n)$.  But directly from the definitions we see that the identity map on elements  gives an isomorphism $\cI_n\cap \cJ_n\iso \cJ_{n-1}$.  So by~\eqref{cJnln} we have
$$
\#(\cI_n\cap \cJ_n)=l_{n-1}
$$
as desired.
\end{proof}

\section{Permutations}
\label{per}

We now prove various log-concavity and log-convexity results concerning subsets of the symmetric group $\fS_n$ of permutations of $[n]$.
Sequences of evaluations of descent and peak polynomials will be shown to be log-concave.  We will also use pattern avoidance to give a third proof of the log-convexity of the Catalan numbers.
In addition to the Order Ideal Lemma, one of our main tools will be the recently-defined middle order distributive lattice on $\fS_n$.

There are two standard partial orders on $\fS_n$: the weak and strong Bruhat orders.
Recently, Bouvel, Ferrari and Tenner~\cite{BFT:bwb} defined a partial order which they call the middle order because it refines the weak order and is refined by the strong.  This order has the advantage of being a distributive lattice and is built using inversions.
Given $\pi=\pi_1\pi_2\ldots\pi_n\in\fS_n$ in one-line notation, its 
{\em set of inversion (values)} is
$$
\Inv\pi=\{(\pi_i,\pi_j) \mid i<j \text{ and } \pi_i>\pi_j\}.
$$
We also say that $\pi_i$ is an {\em inversion top} if $(\pi_i,\pi_j)\in\Inv \pi$ for some $\pi_j$.
For example, if 
\beq
\label{pi:ex}
\pi=415632
\eeq
then
$$
\Inv 415632=\{(4,1),\ (4,2),\ (4,3),\ (5,2),\ (5,3),\ (6,2),\ (6,3),\ (3,2)\}
$$
so $4$ is an inversion top in $3$ inversions, $5$ and $6$ are inversion tops in $2$ inversions each, and $3$ is an inversion top in one inversion.  
Clearly $i$ can be an inversion top in anywhere from $0$ to $i-1$ inversions.
The {\em inversion table} of $\pi$ is
$$
\io(\pi) = (\io_1,\io_2,\ldots,\io_n)
$$
where
$$
\io_i = \text{ \# of inversions in which $i$ is an inversion top}.
$$
Returning to our example
$$
\io(\pi) = (0,0,1,3,2,2).
$$
Let
\beq
\label{cI_n}
\cI_n = \{\io = (\io_1,\io_2,\ldots,\io_n) \mid 
0\le\io_i<i \text{ for all } i\in[n]\}.
\eeq
It is well-known that there is a bijection $\fS_n\ra\cI_n$ given by
$\pi\mapsto \io(\pi)$.

\bfi
\bce
\begin{tikzpicture}
\filldraw(1,0) circle(.1);  
\draw(1,-.5) node{123};
\filldraw(2,1) circle(.1); 
\draw(2.7,1) node{132};
\filldraw(3,2) circle(.1); 
\draw(3.7,2) node{312};
\filldraw(0,1) circle(.1);
\draw(-.7,1) node{213};
\filldraw(1,2) circle(.1); 
\draw(.3,2) node{231};
\filldraw(2,3) circle(.1);
\draw(2,3.5) node{321};
\draw (1,0)--(3,2)--(2,3)--(0,1)--(1,0) (2,1)--(1,2);
\draw(5,1.5) node{$\iso$};
\begin{scope}[shift={(8,0)}]
\filldraw(1,0) circle(.1);  
\draw(1,-.5) node{(0,0,0)};
\filldraw(2,1) circle(.1); 
\draw(3,1) node{(0,0,1)};
\filldraw(3,2) circle(.1); 
\draw(4,2) node{(0,0,2)};
\filldraw(0,1) circle(.1);
\draw(-1,1) node{(0,1,0)};
\filldraw(1,2) circle(.1); 
\draw(0,2) node{(0,1,1)};
\filldraw(2,3) circle(.1);
\draw(2,3.5) node{(0,1,2)};
\draw (1,0)--(3,2)--(2,3)--(0,1)--(1,0) (2,1)--(1,2);
\end{scope}
\end{tikzpicture}
\ece
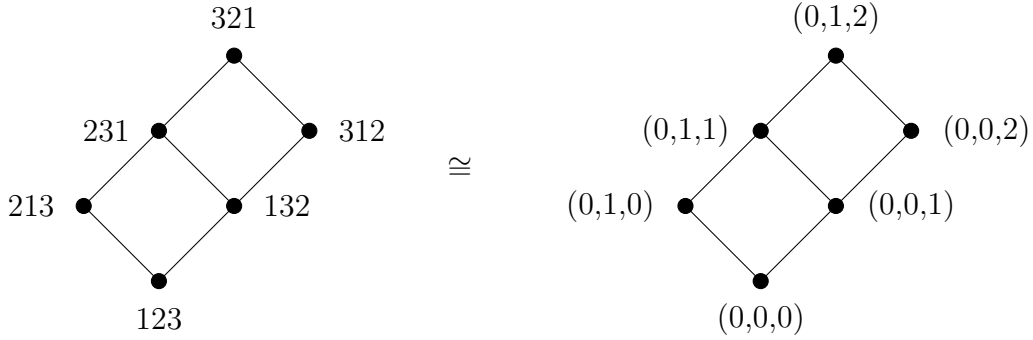
\capt{The middle order on $\fS_3$ both in terms of permutations and inversion tables}
\label{fS_3}
\efi

We can now use the bijection just given to define the {\em middle order}
$(\fS_n,\preceq)$ by
$$
\pi\preceq\si \text{ if and only if } \io(\pi)\le\io(\si)
$$
where the order on inversion tables is component-wise.
Because of~\eqref{cI_n} we have that the middle order is isomorphic to a product of chains
$$
\fS_n \iso [0,0] \times [0,1] \times \cdots \times [0,n-1]
$$
where $[0,i]=\{0,1,\ldots,i\}$ with the usual total order on the integers.
It follows that this order is a distributive lattice.
In Figure~\ref{fS_3} we display this order both on the permutations in $\fS_3$ on the left, as well as on their corresponding inversion tables on the right.
It will sometimes be convenient to work directly with 
$(\cI_n,\le)$ rather than $(\fS_n,\prec)$.  We begin by showing that the sequence of factorials $n!=\#\fS_n$ is log-convex.  Of course, this can be proved by simple arithmetic.  But our proof will serve as a model for later, more complicated ones, using the middle order.
\bth
The sequence $(n!)_{n\ge0}$ is log-convex
\eth
\begin{proof}
Consider the two lower order ideals of $\cI_{n+1}$ given by
$$
I=\{\io\in \cI_{n+1} \mid \io_{n+1} = 0\}
$$
and
$$
J=\{\io\in\cI_{n+1} 
\mid 0\le \io_i < i - 1 \text { for all } 2\le i\le n+1\}.
$$
Now $I\iso\cI_n\iso J$ where the first, respectively second, isomorphism is obtained by removing $\io_{n+1}$, respectively $\io_1$.
Similarly $I\cap J\iso\cI_{n-1}$ and the Order Ideal Lemma once again finishes the demonstration.
\end{proof}

\subsection{Descent polynomials}

We now prove that sequences of evaluations of descent polynomials are log-concave.
After this paper was published, it was discovered that therre was an error in the proof in that the intersection of the two proposed ideals did not have the correct cardinality.  This was corrected in a later corrigendum, and that is the proof which will be presented here.
The {\em descent set} of $\pi\in\fS_n$ is
$$
\Des\pi=\{i \mid \pi_i>\pi_{i+1}\}.
$$
Note that, unlike the inversion set, we are using the positions of the descents.
If $\pi$ is as in~\eqref{pi:ex} then
\beq
\label{Des:ex}
\Des 415632 = \{1,4,5\}
\eeq
Let $S$ be any finite set of positive integers and consider
$$
D_n(S)=\{\pi\in\fS_n \mid \Des\pi = S\}
$$
as well as
$$
d_n(S) = \# D_n(S)
$$
where the latter is called the {\em descent polynomial} corresponding to $S$.
For more information about descent polynomials see~\cite{DLHIOS:dp,ogu:dpp,GG:qdp,ben:scs,ray:gdp,JM:rdp}.
The following is a classic result of MacMahon~\cite{mac:ca}.
\bth[\cite{mac:ca}]
For any set $S$ and all $n>\max S$ we have that $d_n(S)$ is a polynomial in $n$.\hqed
\eth

In order to prove that $(d_n(S))_{n\ge0}$ is always log-concave we will need a variant of the middle order which considers positions rather than values.
The {\em positional inversion table} of $\pi\in\fS_n$ is 
$$
\ka(\pi) = (\ka_1,\ka_2,\ldots,\ka_n)
$$
where
$$
\ka_i = \#\{j>i \mid \pi_j<\pi_i\}.
$$
In other words, $\ka_i$ is the number of inversions with $\pi_i$ as inversion top.  Continuing with our example permutation
\beq
\label{ka:ex}
\ka(415632) = (3,0,2,2,1,0).
\eeq
Clearly $0\le \ka_i \le n-i$ for all $i\in[n]$.  Consider
$$
\cK_n  = \{\ka = (\ka_1,\ka_2,\ldots,\ka_n) \mid 
0\le\ka_i\le n-i \text{ for all } i\in[n]\}.
$$
The bijection $\cK_n\ra\fS_n$ given by $\ka\mapsto\pi$  where $\ka(\pi)=\ka$ will be useful to us, so we describe it explicitly.  Assuming 
that $\pi_1,\ldots,\pi_{i-1}$ have been constructed, we let
\beq
\label{pi:ka}
\pi_i = 
\text{ the $(v_i+1)$st smallest element of $[n]\setm\{\pi_1,\ldots,\pi_{i-1}$\}}.
\eeq
By way of illustration, suppose $\ka$ is as given in~\eqref{ka:ex}.
Since $\ka_1=3$ we let $\pi_1$ be the $(3+1)$st smallest element of $[6]$, that is $\pi_1=4$.  Now $\ka_2=0$ so $\pi_2$ will be the smallest element of $[6]\setm\{4\}$ so that $\pi_2 = 1$.  Next $\ka_3=2$ so we pick the third smallest element of $[6]\setm\{1,4\}$ which gives $\pi_3=5$.  Continuing in this way, we get $\pi=415632$ which agrees with the permutation which started this example.

We now define the {\em $\ka$-middle order} $(\fS_n,\tleq)$ by
$$
\pi\tleq\si \text{ if and only if } \ka(\pi)\le\ka(\si)
$$
component-wise.  In this partial order we have
$$
\fS_n \iso [0,n-1] \times [0,n-2] \times \cdots\times [0,0]
$$
so that, again, we have a distributive lattice.  And, as with middle order, sometimes we will choose to work directly with $\cK_n$ rather than $\fS_n$.
To do this, we need to be able to read off $\Des\pi$ from $\ka(\pi)$ which turns out to be easy to do (and is one of the reasons for using $\ka$ rather than $\io$).  Define the {\em descent set} of $\ka\in\cK_n$ to be
$$
\Des\ka = \{i \mid \ka_i>\ka_{i+1}\}.
$$
Note that we do not include $i$ in the descent set if $\ka_i=\ka_{i+1}$.
If $\ka$ is as in~\eqref{ka:ex} then
$$
\Des\ka=\{1,4,5\}.
$$
Note that this is the same descent set as in~\eqref{Des:ex}.  This is not an accident.
\begin{lem}
\label{Des:ka}
For any $\pi\in\fS_n$ we have
$$
\Des \ka(\pi) = \Des\pi.
$$
\end{lem}
\begin{proof}
Say $\pi=\pi_1\pi_2\ldots \pi_n$ and 
$\ka(\pi)=(\ka_1,\ka_2,\ldots,\ka_n)$.  Suppose first that $i\in\Des \ka(\pi)$ so that $\ka_i>\ka_{i+1}$.  From~\eqref{pi:ka} 
we have that $\pi_i$ is the $(\ka_i+1)$st smallest element of a set $S$  and $\pi_{i+1}$ is the $(\ka_{i+1}+1)$st smallest in $S\setm\{\pi_i\}$.
But since $\ka_{i+1}<\ka_i$ we also have that $\pi_{i+1}$ is the $(\ka_{i+1}+1)$st smallest in $S$ itself.
Using the inequality
$\ka_{i+1}<\ka_i$
again shows that $\pi_{i+1}< \pi_i$ so that $i\in\Des\pi$.  By a similar argument, if $i\not\in\Des \ka(\pi)$
then $i\not\in\Des\pi$ which completes the proof.
\end{proof}

We now define a partial order $(D_n(S),\tleq)$ by restricting the $\ka$-middle order on $\fS_n$ to $D_n(S)$.  We need to show that we still have a distributive lattice.  In fact, we will show that $D_n(S)$ is a sublattice of $\fS_n$ under $\tleq$.  To do so, it will be convenient to extend the minimum and maximum functions to  two integer vectors
$\ka=(\ka_1,\ldots,\ka_n)$ and $\chi=(\chi_1,\ldots,\chi_n)$ component-wise so that 
\beq
\label{min}
\min\{\ka,\chi\} = (\min\{\ka_1,\chi_1\},\ldots, \min\{\ka_n,\chi_n\})
\eeq
and
\beq
\label{max}
\max\{\ka,\chi\} = (\max\{\ka_1,\chi_1\},\ldots, \max\{\ka_n,\chi_n\}).
\eeq
\begin{lem}
If $\ka,\chi\in\cK_n$
have $\Des \ka = S = \Des \chi$ then
$$
\Des(\min\{\ka,\chi\}) = S =\Des(\max\{\ka,\chi\}).
$$
\end{lem}
\begin{proof}
 We will prove the $\min$ statement as the $\max$ one is similar. It suffices to show that if $i\in S$ then $i\in\Des(\min\{\ka,\chi\})$  and similarly for $i\not\in S$.
 
 If $i\in S$ then we have $\ka_i>\ka_{i+1}$ and $\chi_i>\chi_{i+1}$  Without loss of generality, we can assume $\ka_i\le \chi_i$  so that 
 $\min\{\ka_i,\chi_i\} = \ka_i$.  Thus
 $$
 \min\{\ka_{i+1},\chi_{i+1}\}\le \ka_{i+1} < \ka_i = \min\{\ka_i,\chi_i\}
 $$
 which implies that $i\in \Des(\min\{\ka,\chi\})$.  The analogous argument when $i\not\in S$ is left to the reader.
\end{proof}

We now have everything in place to show that $D_n(S)$ is a distributive lattice.
\begin{lem}
For any $S$ the $\ka$-middle order on $D_n(S)$ forms a distributive lattice.    
\end{lem}
\begin{proof}
Consider $\pi,\si\in D_n(S)$.  Define the meet and join operations by letting $\pi\mt\si$ and $\pi\jn\si$ be the unique permutations such  that
$$
\ka(\pi\mt\si) = \min\{\ka(\pi),\ka(\si)\}
$$
and
$$
\ka(\pi\jn\si) = \max\{\ka(\pi),\ka(\si)\}.
$$

We must first check that the meet and join are still in $D_n(S)$, that is, have descent set $S$.  But $\Des\pi=S=\Des\si$ so that, by the previous lemma, 
$\Des \ka(\pi)= S = \Des \ka(\si)$.  Now using both the previous lemmas
$$
\Des \pi\mt\si=\Des \ka(\pi\mt\si)= \Des(\min\{\ka(\pi),\ka(\si)\})=S.
$$
The join case is similar.

We now show that the definition of meet  actually gives a greatest lower bound, leaving  the least upper bound property as an exercise for the reader.
First of all we note that $\ka(\pi\mt\si)= \min\{\ka(\pi),\ka(\si)\}\le \ka(\pi), \ka(\si)$ component-wise so that $\pi\mt\si\le \pi,\si$.
Now suppose that $\tau\le\pi,\si$.  By definition of the $\ka$-middle order we have 
$\ka(\tau)\le \ka(\pi),\ka(\si)$ component-wise.  But then 
$\ka(\tau)\le \min\{\ka(\pi),\ka(\si)\}$ so that $\tau\le \pi\mt\si$.

That meet distributes over join follows easily from the fact that min distributes over max and so the proof is omitted.
\end{proof}

Now define
$$
\cK_n(S) =\{ \ka\in\cK_n \mid \Des\ka=S\}
$$
and partially order this set component-wise.  From what we have shown, $\cK_n(S)\iso D_n(S)$ and so is a distributive lattice of size $d_n(S)$.  We will need the following facts about $\cK_n(S)$.

\ble
Suppose $S\neq\emp$ and let $m=\max S$.  Suppose $\ka=(\ka_1,\ldots,\ka_n)\in\cK_n(S)$.
\ben
\item[(a)]  We have $\ka_m\ge1$.
\item[(b)]  We have $\ka_i=0$ for $i>m$.
\een
\ele
\bprf
(a)  Since $m\in S = \Des\ka$ we must have $\ka_m>\ka_{m+1}\ge0$.  But this implies $\ka_m\ge1$.
\medskip

(b)  Suppose, towards a contradiction, that $\ka_i>0$ for some $i>m$.  By definition of $\cK_n$ we have $\ka_n=0$.  So there must be an index $j\ge i>m$ such that $\ka_j>0$ and $\ka_{j+1}=0$.  But then $j\in\Des\ka$, contradicting the maximality of $m$.
\eprf

We can now prove the main result of this subsection.

\bth
\label{d_n:thm}
For any $S$, the sequence $(d_n(S))_{n\ge0}$ is log-concave.
\eth
\begin{proof}
The result is easy if $S=\emp$, so we assume $S$ is nonempty.
Keeping the notation of the previous lemma, we see that if $\ka\in\cK_n(S)$ then in fact
$$
\ka\in [0, n-1] \times [0, n-2] \times\cdots\times [0, n-m+1] \times [1, n-m] \times [0, 0]^{n-m}
$$
where the superscript denotes the $(n-m)$-fold product.

Now consider the lattice $\cK_{n+1}(S)$ so that $|\cK_{n+1}(S)|=d_{n+1}(S)$.  Inside this lattice, define the set
$$
I=\cK_n(S) \times[0,0].
$$
Clearly $I$ is a lower order ideal with $I\iso D_n(S)$ so that $|I|=d_n(S)$.  Also construct
$$
J = \cK_{n+1}(S) \cap
([1, n] \times [1, n-1] \times \cdots\times [1, n-m+2] \times [2, n-m+1]  \times [0, 0]^{n-m+1}).
$$
From its definition, $J$ is an upper order ideal. To show $|J|=d_n(S)$  we will construct a poset isomorphism $J\iso D_n(S)$.  If $\ka\in  J$ then  factor it as a concatenation $\ka = \la\mu$ where $|\la|=m$ and $|\mu|=n+1-m$.  Map $\ka$ to $\ka' = \la'\mu'$ where $\la'$ is the sequence obtained from $\la$ by subtracting $1$ from each entry, and $\mu'=[0,0]^{n-m}$.  It easy easy to see that this is an isomorphism by constructing the inverse map.

Finally, we have
$$
I\cap J = K_{n+1}(S) \cap
([1, n-1] \times [1, n-2] \times \cdots\times [1, n-m+1] \times [2, n-m]  \times [0, 0]^{n-m+1}).
$$
Using arguments similar to those in the previous paragraph, one shows that $I\cap J\iso D_{n-1}(S)$ which implies $|I\cap J|=d_{n-1}(S)$.  Now we are done by part (b) of the Order Ideal Lemma.
\end{proof}

Note that in the special case $S=[k]$ we have $d(S;n)=\binom{n-1}{k}$ since any $\pi\in D(S;n)$ can be formed by putting $1$ in the $(k+1)$st  position, then choosing $k$ elements of $[2,n]$ to put before $1$ in decreasing order, and arranging the rest after $1$ in increasing order.  Thus we have obtained another proof of Corollary~\ref{BinCo}.

\subsection{Peak polynomials}

We now consider sequences of evaluations of peak polynomials.  The {\em peak set} of
$\pi\in\fS_n$ is
$$
\Pk \pi = \{i \mid \pi_{i-1} < \pi_i > \pi_{i+1}\}.
$$
Using our canonical example~\eqref{pi:ex} we have
$$
\Pk 415632 = \{4\}.
$$
As with the descent set, for a finite set of positive integers $S$ we let
$$
P_n(S) = \{\pi\in\fS_n \mid \Pk\pi = S\}.
$$
There are some sets $S$ such that $P_n(S)=\emp$ for all $n\ge0$.
But if $S$ does not contain $1$ and also does not contain two consecutive integers, one can easily construct $\pi\in\fS_n$ with $\Pk\pi = S$ if $n>\max S$.  We call such $S$ 
{\em admissible}.  The following theorem was proved by Billey, Burdzy, and Sagan~\cite{BBS:pgp}.
\bth
If $S$ is admissible then for $n>\max S$ we have
$$
\#P_n(S) = p_n(S) 2^{n-\#S-1}
$$
where $p_n(S)$ is a polynomial in $n$.\hqed
\eth
The polynomial $p_n(S)$ is called the {\em peak polynomial} and more information about these polynomials can be found in~\cite{BBPS:mew,BFT:crp,DLHIO:ppp,ogu:dpp,GG:qdp}.  We note that the peak set of $\pi$ can be read off from the descent set since
a peak is exactly a non-descent followed by a descent.  And the powers of two in $\#P_n(S)$ will cancel in a log-concavity inequality.  So the proof of the next result closely parallels the one already given for the descent polynomial and  the details are omitted.
\bth
If $S$ is admissible then the sequence $(p_n(S))_{n\ge0}$ is log-concave.\hqed
\eth

\subsection{Pattern avoidance}

In this subsection we will give a third proof that the Catalan sequence is log-convex, this one using pattern avoidance and the  (ordinary) middle order.
Say $\pi\in\fS_n$ contains permutation $\si\in\fS_k$ as a {\em classical pattern} if there is a subsequence $\si'=\pi_{i_1}\pi_{i_2}\ldots \pi_{i_k}$ of $\pi$, called a {\em copy},
whose elements are in the same relative order as those of $\si$.
In a {\em bivincular pattern} certain pairs of elements of $\si'$ are specified to be in adjacent positions (indicated by a vertical bar between the two in $\si$), 
or have adjacent values (indicated by a horizontal bar above the smaller of the two in $\si$). For example, copies of $\si=231$ in $\pi=643512$  are  $451$, $452$, $351$, and $352$.  Only one of these  is a copy of $2|3\ol{1}$ namely $352$.
If $\si$ is any bivincular pattern then its {\em avoidance set} is
$$
\Av_n(\si)=\{\pi\in\fS_n \mid \text{$\pi$ contains no copy of $\si$}\}.
$$
It is well-known that if $\si\in\fS_3$ then
\beq
\label{S_3}
\#\Av_n(\si)=C_n,
\eeq
the $n$th Catalan number.

The following result of   Claesson, Jelínek,  and Steingrímsson~\cite[Lemma 10]{CJS:ubs} will permit us to transfer results from $\fS_n$ to $\cI_n$.
\begin{lem}[\cite{CJS:ubs}]
For all $n\ge0$, we have $\pi\in\Av_n(213)$ if and only if $\io(\pi)$ is weakly increasing.\hqed
\end{lem}

We can now use~\eqref{S_3} to reprove the log-convexity of the Catalan sequence.  
\begin{thm}
The sequence $(\#\Av_n(213))_{n\ge0}=(C_n)_{n\ge0}$ is log-convex.
\end{thm}
\begin{proof}
We first show that the partial order on  $\Av_n(213)$ induced from the middle order is a distributive lattice.  We claim that this order, is in fact, a sublattice of the middle order of $\fS_n$ and thus must be distributive.  The meet and join  of $\pi,\si\in\fS_n$ are obtained by taking pairwise minima and maxima of $\io(\pi)$ and $\io(\si)$.  From this it is easy to see that if $\io(\pi)$ and $\io(\si)$ are weakly increasing then so are $\io(\pi)\mt \io(\si)$ and $\io(\pi)\jn \io(\si)$.  The previous lemma completes the proof of the claim.

We now choose two upper order ideals to finish the demonstration of the theorem.
We use multiplicity notation for inversion sequences in the usual way.
In particular, let
$$
I = \{\io\in \io(\Av_{n+1}(213)) \mid \io \ge (0,1^n)\},
$$
and
$$
J= \{\io\in \io(\Av_{n+1}(213)) \mid \io \ge (0^n,n)\}.
$$
The verification that these ideals and their intersection have the correct sizes is routine and left to the reader.  The Order Ideal Lemma completes the proof.
\end{proof}

We note that one can give a similar proof of the log-convexity of the Catalan sequence using
middle order restricted to $\Av(312)$.  One can not apply the same technique to any other permutation $\pi\in\fS_3$ because, as is clear from Figure~\ref{fS_3}, the restriction of middle order to $\Av(\pi)$ is not a distributive lattice.

\section{Set partitions}
\label{sp}

For our last applications of the Order Ideal Lemma, we use set partitions and noncrossing set partitions to prove log-concavity results about Stirling numbers of the second kind and Narayana numbers.  In both cases, it will be convenient to express the distributive lattices in terms of restricted growth functions.

\subsection{Stirling numbers of the second kind}

We first set notation and basic definitions.  A {\em set partition of $[n]$},
$\be=B_1/B_2/\ldots/B_k\ptn[n]$, is a family of disjoint subsets $B_i$ called {\em blocks} whose disjoint union is $\uplus_i B_i = [n]$.  
In examples, we will eliminate the set braces and commas from the $B_i$.
We will also always  write our partitions in {\em standard form} which means that
$$
1=\min B_1 < \min B_2 <\ldots<\min B_k.
$$
For example, 
\beq
\label{be:ex}
\be= 12359/46/78 \ptn [9].
\eeq
We let
$$
S([n],k) = \{ \be \mid \text{$\be\ptn[n]$ with $k$ blocks}\}.
$$
The {\em Stirling numbers of the second kind} are
$$
S(n,k) =\# S([n],k).
$$

Set partitions are in bijection with certain sequences called restricted growth functions.
A {\em restricted growth function} (RGF) is a sequence $\rho=\rho_1\rho_2\ldots\rho_n$ of positive integers satisfying
\ben
\item $\rho_1 = 1$, and 
\item for $i\ge 2$ we have
$$
\rho_i \le 1 + \max(\rho_1 \rho_2\ldots\rho_{i-1}).
$$
\een
We call $n$ the {\em length of $\rho$} and write $|\rho|=n$
To illustrate
\beq
\label{rho:ex}
\rho = 111212331
\eeq
is an RGF with $|\rho|=9$, while $\tau= 111212431$ is not since $\tau_7=4$ but $1+\max(111212)=3$.
We will use the notation
$$
\RGF(n,k) = \{\rho \mid  |\rho|=n \text{ and } \max\rho = k\}.
$$

\bfi
\bce
\begin{tikzpicture}
\filldraw(2,0) circle(.1);
\draw(2,-.5) node{$1112$};
\filldraw(2,1) circle(.1);
\draw(1.3,1) node{$1121$};
\filldraw(1,2) circle(.1);
\draw(.3,2) node{$1211$};
\filldraw(3,2) circle(.1);
\draw(3.7,2) node{$1122$};
\filldraw(0,3) circle(.1);
\draw(-.7,3) node{$1221$};
\filldraw(2,3) circle(.1);
\draw(2.7,3) node{$1212$};
\filldraw(1,4) circle(.1);
\draw(1,4.5) node{$1222$};
\draw (2,0)--(2,1)--(0,3)--(1,4)--(3,2)--(2,1) (1,2)--(2,3);
\draw(5,2) node{$\iso$};
\begin{scope}[shift={(7,0)}]
\filldraw(2,0) circle(.1);
\draw(2,-.5) node{$(14,11)$};
\filldraw(2,1) circle(.1);
\draw(1,1) node{$(13,11)$};
\filldraw(1,2) circle(.1);
\draw(0,2) node{$(12,11)$};
\filldraw(3,2) circle(.1);
\draw(4,2) node{$(13,12)$};
\filldraw(0,3) circle(.1);
\draw(-1,3) node{$(12,21)$};
\filldraw(2,3) circle(.1);
\draw(3,3) node{$(12,12)$};
\filldraw(1,4) circle(.1);
\draw(1,4.5) node{$(12,22)$};
\draw (2,0)--(2,1)--(0,3)--(1,4)--(3,2)--(2,1) (1,2)--(2,3);    
\end{scope}
\end{tikzpicture}
\ece
\capt{The partial order on $\RGF(4,2)$ both in terms of RGFs and $(F,R)$ pairs}
\label{RGF(4,2)}
\efi

There is a well-known bijection $S([n],k)\ra \RGF(n,k)$ defined by sending
$\be=B_1/B_2/\ldots/B_k$ to $\rho=\rho_1\rho_2\ldots\rho_n$ where
$$
\rho_i = j \text{ if and only if } i\in B_j.
$$
The reader can check that the partition in~\eqref{be:ex} is sent to the RGF in~\eqref{rho:ex} under this map.

To describe the partial order on $\RGF(n,k)$ we will need two sequences.  If $\rho\in\RGF(n,k)$ then its {\em sequence of first occurrences} (firsts) is
$$
F(\rho) = f_1 f_2\ldots f_k
$$
defined by
$$
f_i = \text{$j$ where $\rho_j$ is the first $i$ in $\rho$}.
$$
Note that since $\rho$ is an RGF we always have $1=f_1<f_2<\ldots<f_k$.
We will also use the {\em rest of $\rho$} which is
$$
R(\rho) = \text{ $\rho$  with its first occurrences removed}.
$$
Note that $|R(\rho)| = n-k$. 
Using our example $\rho$, we have
$$
F(111212331) = 147 \text{ and } R(111212331) = 111231.
$$
Finally, we define a partial order  $(\RGF(n,k),\preceq)$ by
$$
\text{$\rho\preceq\tau$ if and only if $F(\rho)\ge F(\tau)$ and $R(\rho)\le R(\tau)$}
$$
where the orders on $F$ and $R$ are component-wise.  Figure~\ref{RGF(4,2)} illustrates this order both on the restricted growth functions  $\rho\in\RGF(4,2)$ on the left and on the pairs $(F(\rho),R(\rho))$ on the right.
\ble
The partial order $(\RGF(n,k),\preceq)$ is a distributive lattice.
\ele
\bprf
We first need to prove the existence of meets and  joins.  As usual, we only do the former as the latter is similar.  Recall the definitions of $\min$ and $\max$ for integer vectors as given in~\eqref{min} and~\eqref{max}.  We will use the notation $F=f_1 f_2\ldots f_k$, 
$R=r_1 r_2\ldots r_{n-k}$, and similarly for other first and rest sequences.

Suppose that $(F,R)$ and $(G,S)$ are two first-rest pairs corresponding to $\rho,\si\in\RGF(n,k)$.  If we can show that 
$(F(\tau),R(\tau))=(\max\{F,G\},\min\{R,S\})$ for some $\tau\in\RGF(n,k)$ then it is easy to see that $\tau$ is the greatest lower bound of $\rho$ and $\si$.
The proof that $\max\{F,G\}$ is still an increasing sequence beginning with $1$ is easy.
So we need to show that if the numbers $1,2,\ldots,k$ are placed in $\min\{R,S\}$ at the indices indicated in $\max\{F,G\}$ to form $\tau$, then their positions will be the first occurrences of these numbers.
This in turn immediately implies that $\tau$ is an RGF.  
Let us consider what happens when inserting $j$.
Without loss of generality, suppose that $f_j\ge g_j$ so that $\max\{f_j,g_j\} = f_j$.  Now consider any $\tau_i$ with $i<f_j$.
There are now two cases depending on whether $\tau_i$ is a first occurrence or not.
If it is, then $\tau_i<j$ as desired because $\max\{F,G\}$ is increasing.
If $\tau_i$ is not a first occurrence then $\tau_i=\min\{r_l,s_l\}$ for some $l\le i$.
So $r_l$ comes before place $f_j$ in $\rho$ which is an RGF, and this implies  $r_l<j$.  Thus
$$
\tau_i=\min\{r_l,s_l\}\le r_l < j
$$
which is again what we wanted.

To finish the proof we note, as usual, that distributivity follows from the properties of $\min$ and $\max$.
\eprf

We now show that a sequence of Stirling numbers of the second kind is log-concave.  Note that this sequence  is different from the one usually studied where $n$ is fixed and $k$ varies which was shown to be log-concave by Lieb~\cite{lie:cpg}.
\bth
For any $k\ge0$, the sequence $(S(n,k))_{n\ge0}$ is log-concave.
\eth
\bprf
Consider $\RGF(n+1,k)$ which, by the previous lemma, is a distributive lattice.  We first look at the subset given by
$$
I=\{\rho\in \RGF(n+1,k) \mid \rho = 11 \rho_3\ldots\rho_{n+1}\}.
$$
We claim that $I$ is a lower order ideal.  For suppose $\rho\in I$ and $\si\preceq\rho$.
But then we have $F(\si)\ge F(\rho) = 1 f_2\ldots f_k$ where $f_2\ge3$ since $\rho$ begins with two $1$'s.
It follows that $F(\si)=1 g_2\ldots g_k$ where $g_2\ge f_2\ge 3$.  So, $\si$ also begins with two
$1$'s as desired.  It should be clear from the definition of an RGF that there is an isomorphism 
$I\iso \RGF(n,k)$ given by removing the first $1$ of $\rho\in I$.

Next consider
$$
J=\{\rho\in \RGF(n+1,k) \mid \rho = \rho_1\rho_2\ldots\rho_n k \text{ and } n+1\not\in F(\rho)\}.
$$
we wish to show that this is an upper order ideal.  So take $\rho\in J$ and $\si\succeq\rho$.
The two conditions on $\rho$ being in $J$ imply that $\rho_{n+1}=k$ is not a first occurrence.  It follows that $R(\rho)= r_1 r_2\ldots r_{n-k} k$.  Now $R(\si)\ge R(\rho)$ implies
that $R(\si)$ also ends with $k$.  And this means that so does $\si$ itself, and that this final $k$ is not a first occurrence since it is from $R(\si)$.  So $\si$ satisfies the two conditions for inclusion in $J$.  There is also an  isomorphism $J\iso \RGF(n,k)$ gotten  by removing the last $k$ of $\rho\in J$ since that $k$ is not a first occurrence.

The reader should now find it easy to prove that $I\cap J \iso \RGF(n-1,k)$.
Thus we are done by the Order Ideal Lemma.
\eprf

\subsubsection{Narayana numbers}

The Narayana numbers can be defined, for $1\le k\le n$, as
$$
N(n,k) = \frac{1}{n}\binom{n}{k-1}\binom{n}{k}.
$$
They refine the Catalan numbers in that
$$
C_n =\sum_{k=1}^n N(n,k)
$$
and count various refinements of the combinatoiral objects enumerated by $C_n$.
We will prove the log-concavity of sequences of Narayana numbers using their interpretation in terms of non-crossing partitions.

\bfi
\bce
\begin{tikzpicture}
\filldraw(2,0) circle(.1);
\draw(2,-.5) node{$1112$};
\filldraw(2,1) circle(.1);
\draw(1.3,1) node{$1121$};
\filldraw(1,2) circle(.1);
\draw(.3,2) node{$1211$};
\filldraw(3,2) circle(.1);
\draw(3.7,2) node{$1122$};
\filldraw(2,3) circle(.1);
\draw(1.3,3) node{$1221$};
\filldraw(2,4) circle(.1);
\draw(2,4.5) node{$1222$};
\draw (2,0)--(2,1)--(1,2)--(2,3)--(3,2)--(2,1) (2,3)--(2,4);    
\draw(5,2) node{$\iso$};
\begin{scope}[shift={(7.5,0)}]
\filldraw(2,0) circle(.1);
\draw(2,-.5) node{$(14,\{\{1,1\}\})$};
\filldraw(2,1) circle(.1);
\draw(.5,1) node{$(13,\{\{1,1\}\})$};
\filldraw(1,2) circle(.1);
\draw(-.5,2) node{$(12,\{\{1,1\}\})$};
\filldraw(3,2) circle(.1);
\draw(4.5,2) node{$(13,\{\{1,2\}\})$};
\filldraw(2,3) circle(.1);
\draw(.5,3) node{$(12,\{\{1,2\}\})$};
\filldraw(2,4) circle(.1);
\draw(2,4.5) node{$(12,\{\{2,2\}\})$};
\draw (2,0)--(2,1)--(1,2)--(2,3)--(3,2)--(2,1) (2,3)--(2,4);    
\end{scope}
\end{tikzpicture}
\ece
\capt{The partial order on $\NC(4,2)$ both in terms of RGFs and $(F,M)$ pairs}
\label{NC(4,2)}
\efi

Call a set partition $\be=B_1/B_2/\ldots/B_k$ {\em crossing} if there exist positive integers $a<b<c<d$ with $a,c\in B_i$ and $b,d\in B_j$ for some $i\neq j$, and {\em non-crossing} otherwise.
Clearly a partition is non-crossing if and only if the associate restricted growth function $\rho=r_1 \ldots r_n$ has no subsequence of the form $ijij$.   
We call such RGFs {\em non-crossing} as well.
For example, in Figure~\ref{RGF(4,2)} on the left, all the partitions are non-crossing except $1212$.
We let
$$
\NC(n,k) = \{\rho\in\RGF(n,k) \mid \text{$\rho$ is non-crossing}\}.
$$
It is well known that
$$
N(n,k) =\#\NC(n,k).
$$

Define $M(\rho)$ to be the multiset underlying $R(\rho)$.  Using the example from the previous subsection
$$
M(111212331) = \{\{ 1^4, 2, 3\}\}.
$$
We now partially order $\NC(n,k)$ by letting $\rho\tleq \si$ if and only if 
\begin{equation}
\label{FM}
F(\rho) \ge F(\si) \text{ and } M(\rho) \le M(\si).
\end{equation}
where we compare two multisets component-wise after writing them out in weakly increasing order.
In Figure~\ref{NC(4,2)} we have written out the order on $\NC(4,2)$ in terms of RGFs (left) and $(F,M)$ pairs (right).

The following lemma will be useful in proving that the partial order on  $\NC(n,k)$ is a lattice.
\begin{lem}
\label{lem:fimi}
Let $F=f_1f_2\ldots f_k$ be an increasing sequence of integers with $f_1=1$ and let 
$M=\{\{ 1^{m_1}, 2^{m_2},\ldots,k^{m_k}\}\}$ be a multiset.  There is a an RGF $\rho$ with $F(\rho)=F$ and $M(\rho)=M$ if and only if for all $i$ with $1<i\le k$ we have
\begin{equation}
\label{fimi}  
f_i-i \le m_1 + m_2 + \cdots + m_{i-1}.
\end{equation}
Furthermore, there is a unique such non-crossing $\rho$.
\end{lem}
\begin{proof}
Suppose first that $\rho$ exists.  Then $f_i-i$ is the number places in $\rho$ before $f_i$ which do not contain a first element.  And since $\rho$ is an RGF, these places can only contain elements smaller than $i$.  Since the sum is the total number of elements smaller than $i$, the inequality follows.  

Now suppose the inequalities hold.  We will build the desired RGF $\rho$.
An example of this construction follows the proof.
First put $1,2,\ldots,k$ in the spaces dictated by $F$.  Now fill each remaining space from left to right while removing elements from $M$ as follows.  To fill a given space, find the nearest first occurrence from $F$ to its left, say that is $i$.  Now pick the largest element of the current version of $M$ which is at most $i$, say that is is $j\le i$.
Finally put a $j$ into the space of $\rho$ and remove a copy of $j$ from $M$.
We see from the arguments in the previous paragraph that there will always be an element available to fill the space and so this algorithm terminates with all spaces filled.  Furthermore, the RGF restriction is satisfied by construction.

We claim that the RGF $\rho$ constructed in the previous paragraph is noncrossing.  Suppose, to the contrary, the it contains a copy of $ijij$.  There are now two cases depending on the relative size of $i$ and $j$.  But they are similar so we will only do the case $i<j$.  Note that the second $i$ and $j$ in the copy are not first occurrences in $\rho$ and so were added during the procedure using $M$.  Let $k$ be the element in a position of $F$ which is closest on the left to the second $i$ in $ijij$.  So $i$ was the largest element of the current $M$ smaller than or equal to $k$.  But $j>i$ and, by definition of an RGF and the choice of $k$, all elements up to the $i$ in question are at most $k$.
In particular $j\le k$ because of the first $j$ in $ijij$.  It follows that there was no $j$ in $M$ when the second $i$ was chosen since the algorithm always picks the largest possible element.  But now it is not possible to pick an element in a later place for the second $j$.

For uniqueness, suppose that there is another associated noncrossing RGF $\rho'$.
Consider the first place where $\rho$ and $\rho'$ differ and let $i'$ be the element of $\rho'$ in that place. The elements of $F$ are the same in both, so $i'$ must be a place filled by $M$ in $\rho$.
Let $k$ be the first occurrence closest to $i'$ on its left.
By definition of how elements are chosen in $\rho$ and the fact that this is the first place where the two RGFs differ, there must be $j'\in M$ with $i'<j'\le k$ which comes later in $\rho'$.  But from what we have established, neither the occurrence of $i'$, nor that of $j'$ is first.
It is now easy to see that $\rho'$ has a copy of $i'j'i'j'$ where the first $i'$, $j'$ are in positions indexed by $F$ and the second two are as constructed.
\end{proof}

To illustrate the building of the RGF  $\rho$ in the previous proof, suppose that
$$
(F,M) = (137,\ \{\{1^3,2^2,3\}\})
$$
We are to build $\rho$ with $|\rho|=|F|+|M|= 3+6 = 9$.  We start by putting $1,2,3$ in the spaces dictated by $F$ to get
$$
\rho= 1\ \_\ 2\ \_\ \_\ \_\ 3\ \_\ \_.
$$
To fill the first space, we look for the closest first occurrence to its left which is $1$.
So  we must choose the largest element of $M$ which is at most $1$.  Of course, this is $1$ itself, so moving a $1$ from $M$ to $\rho$ gives
$$
\rho= 1\ 1\ 2\ \_\ \_\ \_\ 3\ \_\ \_ \text{ and } M=\{\{1^2,2^2,3\}\}.
$$
The next space is closest on the left to first occurrence $2$.  Since $M$ contains a $2$, we move it to $\rho$ to obtain
$$
\rho= 1\ 1\ 2\ 2\ \_\ \_\ 3\ \_\ \_ \text{ and } M=\{\{1^2,2,3\}\}.
$$
Similarly the next space is filled with a $2$
$$
\rho= 1\ 1\ 2\ 2\ 2\ \_\ 3\ \_\ \_ \text{ and } M=\{\{1^2,3\}\}.
$$
Now $M$ no longer contains any $2$'s so the next space gets the largest element of $M$ less than $2$ which is a $1$
$$
\rho= 1\ 1\ 2\ 2\ 2\ 1\ 3\ \_\ \_ \text{ and } M=\{\{1,3\}\}.
$$
The final two spaces have $3$ as their closest first occurrence, and both elements of $M$ are at most $3$.  So they are added to $\rho$ in decreasing order to finally obtain
$$
\rho= 1\ 1\ 2\ 2\ 2\ 1\ 3\ 3\ 1.
$$

The next lemma will help us show that meets and joins exist in the partial order on $\NC(n,k)$.
\begin{lem}
Let $(F,M)$ and $(F',M')$ be set-multiset pairs satisfiying equation~\eqref{fimi}.
Then the pairs $(\max\{F,F'\},\min\{M,M'\})$ and
$(\min\{F,F'\},\max\{M,M'\})$ satisfy the same inequalities.
\end{lem}
\begin{proof}
We will only provide the details for $(\max\{F,F'\},\min\{M,M'\})$ as  the other case is similar.  Let $F'' = \max\{F,F'\}$ and $M''=\min\{M,M'\}$.  As usual, we use the notation 
$F=f_1f_2\ldots f_k$, $M=\{\{1^{m_1}, 2^{m_2},\ldots, k^{m_k}\}\}$, and similarly with primes or double primes for the other pairs.  Without loss of generality we can assume that $f_i\ge f'_i$ so that
$f''_i=f_i$.  Also $M\ge M''$ component wise so that the number of elements less than $i$ in $M$ must be at at most their number in $M''$, that is,
$$
m_1 + m_2 + \cdots + m_{i-1} \le m''_1 + m''_2 +\cdots + m''_{i-1}.
$$
Thus
$$
f''_i -i = f_i -i \le m_1 + m_2 + \cdots + m_{i-1} \le m''_1 + m''_2 +\cdots + m''_{i-1}
$$
as desired.
\end{proof}

\begin{lem}
The poset $(\NC(n,k),\tleq)$  is a distributive lattice.    
\end{lem}
\begin{proof}
Suppose $\rho,\si\in\NC(n,k)$.  By Lemma~\ref{lem:fimi}, the pairs $(F(\rho),M(\rho))$ and $(F(\si),M(\si))$ satisfy the inequalities~\eqref{fimi}.
And by the previous lemma, \eqref{fimi} is still satisfied by the pair
$(\max\{(F\rho),F(\si)\},\ \min\{M(\rho),M(\si)\})$.  Using Lemma~\ref{lem:fimi} again, there is a unique noncrossing RGF $\rho\mt\si$ with 
$$
F(\rho\mt\si) =  \max\{F(\rho),F(\si)\} \text{ and } M(\rho\mt\si) =  \min\{M(\rho),M(\si)\}.
$$

We claim that $\rho\mt\si$ is the greatest lower bound of $\rho$ and $\si$.
Indeed, 
$F(\rho\mt\si)= \max\{F(\rho),F(\si)\}\ge F(\rho),F(\si)$
 and  $M(\rho\mt\si)= \min\{M(\rho),M(\si)\}\le M(\rho),M(\si)$
so that
$\rho\mt\si\tleq \rho,\si$ by~\eqref{FM}.
And if $\tau\tleq\rho,\si$ then  $F(\tau)\ge F(\rho),F(\si)$ and $M(\tau)\le M(\rho),M(\si)$
by~\eqref{FM} again.  So $F(\tau)\ge \max\{F(\rho),F(\si)\}$ and $M(\tau)\le\min\{M(\rho),M(\si)\}$.
Using the definition of the partial order one last time we get $\tau\tleq\rho\mt\si$ as desired.

Similarly, one constructs $\rho\jn\si$ as the unique noncrossing RGF with
$$
F(\rho\jn\si) =  \min\{F(\rho),F(\si)\} \text{ and } M(\rho\jn\si) =  \max\{M(\rho),M(\si)\}
$$
and shows that it is indeed a least upper bound.  The fact that the resulting lattice is distributive follows from the fact that max distributes over min and vice-versa.
\end{proof}

We can finally prove our main result of this subsection.
\begin{thm}
For fixed $k$, the sequence of Narayana numbers $(N(n,k))_{n\ge0}$ is log-concave.  
\end{thm}
\begin{proof}
Throughout this proof we will consider the elements $\rho\in\NC(n,k)$ as set-multiset ordered pairs $(F(\rho),M(\rho))$ where $F(\rho)=f_1f_2\ldots f_k$ and $M(\rho)=\{\{1^{m_1},2^{m_2},\ldots,k^{m_k}\}\}$.  

Consider the  subset of $\NC(n+1,k)$ given by
$$
I = \{ (F,M) \in \NC(n+1,k) \mid f_2\ge 3\}.
$$
Note that $I$ is a lower order ideal since if $(F,M)\in I$ and $(F',M')\le (F,M)$ then $F'\ge F$.  It follows that
$f_2'\ge f_3\ge 3$ so that $(F',M')\in I$.  We also claim that $I\iso\NC(n,k)$.  If $(F,M)\in I$ then $f_2\ge3$.  
Thus $\rho$, the associated RGF, must begin with at least two $1$'s.  This forces $m_1\ge 1$.
Now map $(F,M)$ to the pair $(F',M')$ where
$$
f_i' = \case{1}{if $i=1$,}{f_i - 1}{if $i\ge2$,}
$$
and
$$
m_i' = \case{m_1-1}{if $i=1$,}{m_i}{if $i\ge2$.}
$$
It is easy to check that this gives an isomorphism.  So we have $\#I = N(n,k)$.

Now define
$$
J  = \{ (F,M) \in \NC(n+1,k) \mid m_k\ge 1\}.
$$
This forces $J$ to be an upper order ideal, for suppose $(F,M)\in J$ and $(F',M')\ge (F,M)$.
Since $m_k\ge1$, the nondecreasing rearrangement of $M$ must end in $k$.  So $M'\ge M$ means the same is true for $M'$, and $m_k'\ge1$.  Again, there is an isomorphism between $J$ and $\NC(n,k)$.
The fact that  $m_k\ge 1$ implies $f_k\le n$ since there are at least two $k$'s in the associated $\rho$ and so the first one can not be in the final position $n+1$.  So the following function is well defined.
Map $(F,M)$ to $(F',M')$ where $F'=F$ and 
$$
m_i' = \case{m_k-1}{if $i=k$,}{m_i}{if $i<k$.}
$$
Again, it is not hard to show that this is an isomorphism.  
The fact that $I\cap J \iso \NC(n-1,k)$  is similarly left to the reader.
\end{proof}

\section{Future directions}
\label{fd}

The purpose of this section is threefold.  First, we will give a proof of the Order Ideal Lemma using the FKG inequality.  Next, we will describe various sequences to which it might be possible to apply our method but which have so far resisted proof.  Finally, we end with some possible avenues for extending the Order Ideal Lemma.

\subsection{Proof of the Order Ideal Lemma}
\label{poi}

Let $\bbR_{\geq0}$ denote the nonnegative real numbers.  For any poset $(P,\preceq)$ we call a  function $f:P\ra\bbR_{\geq0}$ {\em increasing} if
$$
x\preceq y \text{ implies }f(x)\le f(y).
$$
Similarly define $f$ to be {\em decreasing}.
Now suppose $(L,\preceq)$ is a lattice.  A function $\mu:L\to \mathbb{R}_{\geq 0}$ is {\em log-supermodular} if for all $x,y\in L$ we have
$$
\mu(x)\mu(y)\leq \mu(x\land y)\mu(x\lor y).
$$
Finally, given a  function $f:L\ra\bbR_{\geq0}$ and a log-supermodular function
$\mu:L\to \mathbb{R}_{\geq 0}$,
we use the notation
$$
S(f) = \sum_{x\in L} f(x)\mu(x).
$$
In particular, if $g:L\ra\bbR_{\geq0}$ is another function then
$$
S(fg) = \sum_{x\in L} f(x)g(x)\mu(x),
$$
and, letting ${\bf 1}:L\ra \bbR_{\geq0}$ be the function defined by 
${\bf 1}(x) = 1$ for all $x\in L$,
$$
S({\bf 1})= \sum_{x\in L} \mu(x).
$$
We can now state the FKG inequality~\cite{FKG:cip}.

\begin{thm}[\cite{FKG:cip}]
\label{FKG}
Suppose $(L,\preceq)$ is a finite distributive lattice and 
$\mu:L\to\mathbb{R}_{\geq 0}$ is a log-supermodular function. Suppose also that 
$f,g: L\to \bbR_{\geq0}$ are two functions.
\ben
\item[(a)]  If $f,g$ are both increasing or both decreasing then
$$
S(f) \cdot S(g) \le S(fg) \cdot S({\bf 1}).
$$
\item[(b)]  If one of $f,g$ is increasing and the other decreasing then

\vs{10pt}

\eqqed{S(f) \cdot S(g) \ge S(fg) \cdot S({\bf 1}).}

\een
\end{thm}

The Order Ideal Lemma, which is restated here for convenience, is a corollary.
\ble[The Order  Ideal Lemma]
Let $L$ be a distributive lattice and suppose that $I,J\sbe L$ are ideals.
\ben
\item[(a)]  If $I,J$ are both lower ideals or both upper ideals then
$$
|I| \cdot |J| \le |I\cap J| \cdot |L|.
$$
\item[(b)]
If one of $I,J$ is a lower order ideal and the other is upper then
$$
|I| \cdot |J| \ge |I\cap J| \cdot |L|.
$$
\een
\ele
\bprf
We will prove (b) as (a) is similar.  Suppose that $I$ is the lower order ideal and $J$ is the upper.  Define indicator functions $f,g: L\to \bbR_{\geq0}$ by
$$
f(x) = \case{1}{if $x\in I$,}{0}{else,}
$$
and
$$
g(x) = \case{1}{if $x\in J$,}{0}{else.}
$$
Since $I$ is a lower order ideal we have that $f$ is decreasing and, by the same token, $g$ is increasing.  Also define  $\mu={\bf 1}$ which is clearly log-supermodular.  Now, by the FKG Theorem, part (b),
$$
|I|\cdot |J| = S(f) \cdot S(g) \ge S(fg) \cdot S({\bf 1}) = |I\cap J| \cdot |L|
$$
which finishes the proof.
\eprf

\subsection{Other sequences}

In this subsection we discuss various sequences to which we hope the Order Ideal Lemma might be applied.

A superset of the set of permutations (viewed as permutations matrices) is the set of alternating sign matrices.  An {\em alternating sign matrix} or ASM is a matrix  such that 
\ben
\item every entry  is $\pm 1$ or $0$, and
\item in each row and each column the nonzero entries alternate $1$ and $-1$, beginning and ending with $1$.
\een
The number of $n\times n$ ASM's is given by 
$$
\asm_n = \prod_{i=0}^{n-1} \frac{(3i+1)!}{(n+i)!}.
$$
This formula was conjectured by Mills, Robbins, and Rumsey in 1983~\cite{MRR:asm} and then given two proofs in 1996, the first by Zeilberger~\cite{zei:pas} and the second by Kuperberg~\cite{kup:apa}.  The same numbers count descending plane partitions, totally symmetric self-complementary plane partitions, configurations of square ice with domain wall boundary conditions, and various other mathematical objects.  Using cancellation of factorials, it is easy to prove the following result.
\bth
\label{asm:thm}
The sequence $(\asm_n)_{n\ge0}$ is log-convex.\hqed
\eth

It is a natural question to find a proof of the previous resul using the Order Ideal Lemma and this has been done by Fischer, Konvalinka, Sagan, and Weigandt~\cite{FKSW:lcc} using the corresponding corner sum matrices.  These ideas can also be used to prove log-concavity and log-convexity for related objects which do not have a nice product formula.

A famous set of permutations of a multiset is the set of parking functions.  A {\em parking function} of length $n$ is a sequence of  integers $\phi=\phi_1\phi_2\ldots\phi_n$ whose weakly increasing rearrangement $\psi=\psi_1\psi_2\ldots\psi_n$ satisfies
$$
1\le \psi_i \le i
$$
for $i\in[n]$.  Otherwise put, and using the notation of~\eqref{la+n} and~\eqref{cI_n}, 
$\psi-1\in\cI_n$.  Parking functions were first defined by Konheim and Weiss~\cite{KW:oda} and have since been widely studied, in part because of their connections with noncrossing partitions, hyperplane arrangements, and other combinatorial constructs.  The number of parking functions of length $n$ is given by
$$
\prf_n = (n+1)^{n-1}.
$$
We note that this is also the number of labeled trees on $n+1$ vertices.
The next result follows from an algebraicly derived theorem of Chen, Wang, and Yang~\cite{CWY:rrs} about strongly $q$-log-concave sequences (see the next section for the definition) applied to certain rooted trees.
\bth[\cite{CWY:rrs}]
The sequence $(n^{n-1})_{n\ge 2}$ is log-convex.\hqed
\eth

Using simple division and the result just given, we obtain the following.
\bth[\cite{CWY:rrs}]
The sequence $(\prf_n)_{n\ge 1}$ is log-convex.\hqed
\eth
It would be very interesting to give a combinatorial proof of this theorem using the Order Ideal Lemma.  Unfortunately, we have been unable to find the appropriate distributive lattice.

Given a permutation, one can look at the values of the peaks rather than their indices.  Specifically, the {\em pinnacle set} of $\pi\in\fS_n$ is
$$
\Pin\pi=\{\pi_i \mid \pi_{i-1}<\pi_i>\pi_{i_1}\}.
$$
Pinnacle sets were first studied by Strehl~\cite{str:eap} although he called them peak sets.  They were then rediscovered by Davis, Nelson, Petersen, and Tenner~\cite{DNPT:psp} and have since been studied by a number of authors~\cite{rus:spf,DHHIN:fep,RT:apo,fang:ere,DLMSSS:psp,min:frp,FNT:psr}.  Interestingly, they do not behave like peaks.  For example, there does not seem to be any corresponding pinnacle polynomial.

Since pinnacles involve permutation values, it is natural to try and study them using the ordinary middle order.  However, the posets obtained by restricting this order to all permutations with a given length and pinnacle set do not always have connected Hasse diagrams.  Instead, consider the following variant which has occurred, for example, in the work of Rusu and Tenner~\cite{RT:apo}.  Let $\si$ be a permutation of distinct positive integers.  Let
$$
\Pin_n(\si) = \{\pi\in\fS_n \mid
\text{the pinnacles of $\pi$ are the elements of $\si$ in that order}\}.
$$
and
$$
\pin_n(\si) = \#\Pin_n(\si).
$$
Computation of examples raise the following questions.
\begin{ques}
For all $n$ and $\si$:
\ben
\item[(a)]  Is the restriction of the middle order to $\Pin_n(\si)$ a distributive lattice?
\item[(b)] Is the sequence $(\pin_n(\si))$ log-concave?
\een
\end{ques}

It would be interesting if one could apply our methods to prove log-convexity of  other pattern avoidance sequences.  In fact, they are conjectured to satisfy a stronger property.
Call $(a_n)_{n\ge0}$ a {\em Stieltjes moment sequence}  if there exists a nonnegative measure $\mu$  supported on $[0,\infty)$  such that
$$
a_n = \int_0^\infty x^n d\mu(x)
$$
for all $n\ge0$.  This is equivalent to requiring that all of the determinants of the Hankel matrix
$$
H= [a_{i+j}]_{i,j\in\bbN}=
\left[
\barr{cccc}
a_0     & a_1   & a_2   & \cdots\\
a_1      & a_2   & a_3   & \cdots\\
a_2      & a_3     & a_4   & \cdots\\
\vdots  &\vdots &\vdots &\ddots
\earr
\right]
$$
are non-negative.  Note that this immediately implies log-convexity by taking adjacent $2\times 2$ minors.  In his thesis, Elvey Price posed the following question.
\begin{ques}[\cite{elv:spe}]
Suppose  $\si\in\fS_n$  for $n\ge2$.  Is the sequence $(\#\Av_n(\si))_{n\ge0}$ always a  Stieltjes moment sequence?
\end{ques}

It is known that the answer to this question is ``yes" for all $\si\in\fS_n$ where $2\le n\le 4$ except possibly $\si=1324$ and the other $\si$ giving the same avoidance sequence.  See the article of Vatter~\cite{vat:app} for details as well as more intriguing open questions related to pattern avoidance.

Finally, given the nice behaviour of the Stirling numbers of the second kind, one could ask what happens with those of the first.  Recall that the {\em signless Stirling numbers of the first kind} are
$$
c(n,k)=\#\{\pi\in\fS_n \mid
\text{$\pi$ has $k$  cycles in its disjoint cycle decomposition}\}.
$$
We have checked the following conjecture for $1\le k\le n\le 100$.
\bcon
Given $k$, there is an integer $N_k$ such that $(c(n,k))_{n\ge0}$ is log-concave for $n<N_k$ and log-convex for $n\ge N_k$.
\econ

\subsection{$q$-analogues, $\bx$-analogues, and total positivity}
\label{qa}

Let $q$ be a variable.
There is a $q$-analogue of the FKG inequality which makes it possible to prove generalizations of our results in a straightforward manner.  Put a partial order $\leq$ on the polynomial algebra $\bbR[q]$ by defining
\beq
\label{q-ord}
p(q)\le s(q) \text{ if and only if } s(q)-p(q)\in\bbR_{\ge0}[q].
\eeq
We can now define a sequence of such polynomials $(p_n(q))_{n\ge0}$ to be {\em $q$-log-concave} if
$$
p_n(q)^2 \ge p_{n-1}(q) p_{n+1}(q)
$$
for all $n\ge1$ and similarly for {\em $q$-log-convexity}.

Suppose   $(L,\preceq)$ is a lattice and $\bbN$ is the nonnegative integers.  Call a function $r:L\ra\bbN$ {\em modular} if
$$
r(x)+r(y) = r(x\mt y) + f(x\jn y).
$$
If $L$ is distributive, then it has a rank function which is an example of a modular function.  If, in addition to $r$, we have a function  $f:L\ra\bbR_{\ge0}$ and a log-supermodular function $\mu:L\ra\bbR_{\ge0}$ then we have an associated polynomial
$$
S[f] = \sum_{x\in L} f(x) \mu(x) q^{r(x)}.
$$
The following result was proved by Bj\"orner~\cite{bjo:qFKG} for the rank function. 
 Then
Chen, Pak, and Panova~\cite{CPP:epi} noted that the same proof would work for any modular function.
\begin{thm}[\cite{bjo:qFKG,CPP:epi}]
\label{qFKG}
Suppose $(L,\preceq)$ is a finite distributive lattice, 
$\mu:L\to\mathbb{R}_{\geq 0}$ is a log-supermodular function,
and $r:L\ra\bbN$ is modular.
Suppose also that 
$f,g: L\to \bbR_{\geq0}$ are two functions.
\ben
\item[(a)]  If $f,g$ are both increasing or both decreasing then
$$
S[f] \cdot S[g] \le S[fg] \cdot S[{\bf 1}].
$$
\item[(b)]  If one of $f,g$ is increasing and the other decreasing then

\vs{10pt}

\eqqed{S[f] \cdot S[g] \ge S[fg] \cdot S[{\bf 1}].}

\een
\end{thm}

The following $q$-analogue of the Order Ideal Lemma follows from the previous theorem in much the same way that the original result follows from the FKG inequality.  So the demonstration is omitted.
For a modular function $r$ will use the notation
$$
[I]_q = \sum_{x\in I} q^{r(x)}.
$$
\ble[The $q$-Order Ideal Lemma]
Let $L$ be a distributive lattice, $r:L\ra\bbN$ be modular,  and suppose that $I,J\sbe L$ are ideals.
\ben
\item[(a)]  If $I,J$ are both lower ideals or both upper ideals then
$$
[I]_q \cdot [J]_q \le [I\cap J]_q \cdot [L]_q.
$$
\item[(b)]
If one of $I,J$ is a lower order ideal and the other is upper then

\vs{10pt}

\eqqed{
[I]_q \cdot [J]_q \ge [I\cap J]_q \cdot [L]_q.
}
\een
\ele

We now show how to use the $q$-Order Ideal Lemma to obtain a generalization of Theorem~\ref{sla} to principal specializations of Schur functions, noting that a similar extension works for the more general Theorem~\ref{Om:thm}.
\bth
For any partition $\la$ the sequence $(s_\la(q,q^2,\ldots,q^n))_{n\ge0}$ 
is $q$-log-concave.
\eth
\bprf
For an SSYT $T$ of shape $\la$ we let
$$
r(T) = \sum_{(i,j)\in\la} T_{i,j}.
$$
Note that this is a modular function since the meet and join of tableaux  are obtained by taking element-wise minima and maxima, and for any real numbers $x,y$ we
have
$$
x+y = \min\{x,y\} + \max\{x,y\}
$$
Letting
$$
\SSYT_\la(n) = \{T\in\SSYT_\la \mid \max T \le n\}
$$
we have
$$
s_\la(q,q^2,\ldots,q^n)
=
\sum_{T\in\SSYT_\la(n)}\ \prod_{(i,j)\in\la} q^{T_{i,j}}
= \sum_{T\in\SSYT_\la(n)} q^{r(T)}.
$$
Now using using the lower order ideal
$$
I = \{ T\in\SSYT_\la(n+1) \mid \max T \le n\},
$$
the upper order ideal
$$
J = \{ T\in\SSYT_\la(n+1) \mid \min T \ge 2\},
$$
and the $q$-Order Ideal Lemma completes the proof.
\eprf

There is an extension of $q$-log-concavity to several variables which would also be worth considering.  Suppose that we have a sequence of positive real numbers $(a_n)_{n\ge0}$ which is log concave.  Rewriting the inequalities gives
$$
\frac{a_0}{a_1} \le \frac{a_1}{a_2} \le \frac{a_2}{a_3} \le \ldots\ .
$$
Now cross-multiplying any two fractions gives the seemingly stronger, but actually equivalent, condition that
$$
a_m a_n \ge a_{m-1} a_{n+1}
$$
for all $0<m\le n$.
We note that the analogous inequalities for polynomials in $q$ are not equivalent to $q$-log-concavity and are  called {\em strong $q$-log-concavity}.

Suppose that $\leq$ is a partial order on  the polynomial ring $\bbR[\bx]$ where  $\bx=\{x_1,x_2,\ldots\}$.  Say that this order is {\em standard} if it satisfies the following three axioms.
\ben
\item If $f(\bx)\in\bbR_{\geq0}[\bx]$ then $f(\bx)\ge0$.
\item If $f(\bx)\le g(\bx)$ then $f(\bx)+ h(\bx)\le g(\bx)+h(\bx)$ for all 
$h(\bx)\in\bbR[\bx]$.
\item If $f(\bx)\le g(\bx)$ then $f(\bx) h(\bx)\le g(\bx)h(\bx)$ for all 
$h(\bx)\in\bbR_{\geq0}[\bx]$.
\een
Note that the order defined by~\eqref{q-ord} with $q$ replaced by $\bx$ everywhere is standard.  Call a sequence
$(f_n(\bx))_{n\ge0}$ {\em strongly $\bx$-log-concave} with respect to a standard partial order if, for all $0<m\le n$, we have 
$$
f_m(\bx) f_n(\bx) \ge f_{m-1}(\bx) f_{n+1}(\bx).
$$
Sagan~\cite{sag:lcs} proved  the following result.
\bth[\cite{sag:lcs}]  Fix $k\ge0$ and let $\le$ be a standard partial order on $\bbR[\bx]$. 
If the sequence $(x_n)_{n\geq0}$ is strongly $\bx$-log-concave then so are the sequences
$$
(e_k(x_1,x_2,\ldots,x_n))_{n\geq0}
\text{ and }
(h_k(x_1,x_2,\ldots,x_n))_{n\geq0}
$$
where $e_k$ and $h_k$ are the $k$th elementary and complete homogeneous symmetric functions, respectively.\hqed
\eth

We note that Chan and Pak~\cite[Theorem 6.1]{CP:mci} have proved a multivariate version of the Ahlswede Daykin (AD) inequality~\cite{AD:iwt}.  Since the ordinary AD inequality can be used to prove FKG, the Chan-Pak result implies  a multivariate analogue of the Order Ideal Lemma which could yield interesting results.

It would be very interesting to prove an $\bx$-analogue of the FKG inequality and apply it to prove (strong) $\bx$-log-concavity and $\bx$-log-convexity results.

Another way to generalized our results would be in terms of total positivity.  Note that we can rewrite the log-concavity of a real sequence $(a_n)_{n\ge0}$  in terms of a determinant
$$
\left|
\barr{cc}
a_n & a_{n+1}\\
a_{n-1} & a_n
\earr
\right|\ge0.
$$
Letting $a_n=0$ for $n<0$, consider the infinite Toeplitz matrix
$$
T= [a_{j-i}]_{i,j\in\bbN}=
\left[
\barr{cccc}
a_0     & a_1   & a_2   & \cdots\\
0       & a_0   & a_1   & \cdots\\
0       & 0     & a_0   & \cdots\\
\vdots  &\vdots &\vdots &\ddots
\earr
\right].
$$
Call $(a_n)_{n\ge0}$ {\em totally positive} if all minors of $T$ are nonnegative.
Note that this implies log-concavity because of the connected $2\times 2$ minors (that is, those consisting of two adjacent rows and two adjacent columns).  It would be very interesting to find an extension of the FKG inequality and the Order Ideal Lemma which could be used to prove total positivity.

\vs{30pt}

{\em Acknowledgment.}  We would like to thank Benjamin Adenbaum and Thomas McConville for useful discussions.  In particular, Adenbum gave references concerning when the weak and strong Bruhat order restricted to an avoidance class is a lattice.  We are also indebted to Vincent Vatter for pointing out the connections with Stieltjes moment sequences.   
We thank Nick Ovenhouse and Wonmyeong Cho for pointing out that the proof of Theorem~\ref{d_n:thm} needed fixing.
Finally, the referee provided useful suggestions and references.

\vs{30pt}

{\em Funding Declaration.}  Neither Jinting Liang nor Bruce Sagan received funding for the writing of this article.

\nocite{*}
\bibliographystyle{alpha}

\newcommand{\etalchar}[1]{$^{#1}$}

\end{document}